\documentclass[11pt]{article}

\usepackage{amstext,amssymb,amsmath,amsbsy}
\usepackage{hyperref}
\usepackage{amscd}
\usepackage{amsfonts}
\usepackage{indentfirst}
\usepackage{verbatim}
\usepackage{amsmath}
\usepackage{amsthm}
\usepackage{enumerate}
\usepackage{graphicx}
\usepackage{color}
\usepackage[OT1]{fontenc}
\usepackage[latin1]{inputenc}
\usepackage[english]{babel}
\usepackage{amssymb}
\usepackage{subfig}
\usepackage{algorithm}
\usepackage{algpseudocode}
\usepackage[shortlabels]{enumitem}

\newcommand{\R}{\mathbb{R}}
\newcommand{\ds}{\displaystyle}

\newcommand{\x}{{\bf x}}
\newcommand{\p}{{\bf p}}

\newcommand{\bn}{{\bf n}}
\newcommand{\bm}{{\bf m}}
\newcommand{\bN}{{\bf N}}

\newtheorem{Remark}{Remark}[section]

\newtheorem*{Assumption*}{Assumption}

\newtheorem{Problem}{Problem}[section]
\newtheorem*{Problem*}{Problem}
\numberwithin{equation}{section}

\usepackage{matlab-prettifier}

\begin{document}

\title{The dimensional reduction method for solving a nonlinear inverse heat conduction problem with limited boundary data}

\author{
Dinh-Nho H\`ao\thanks{%
Institute of  Mathematics, Vietnam Academy of Science and Technology, \texttt{hao@math.ac.vn}.}
\and Thuy T. Le\thanks{%
Department of Mathematics and Statistics, University of North Carolina at
Charlotte, Charlotte, NC 28223, USA, \texttt{tle55@uncc.edu}, corresponding author} 
 \and Loc H. Nguyen\thanks{%
Department of Mathematics and Statistics, University of North Carolina at Charlotte, Charlotte, NC 28223, USA, \texttt{lnguye50@charlotte.edu}.} }
\date{}
\maketitle
\begin{abstract}
	
The objective of this article is to introduce a novel technique for computing numerical solutions to the nonlinear inverse heat conduction problem. This involves solving nonlinear parabolic equations with Cauchy data provided on one side $\Gamma$ of the boundary of the computational domain $\Omega$.
The key step of our proposed method is the truncation of the Fourier series of the solution to the governing equation. 
The truncation technique enables us to derive a system of 1D ordinary differential equations.
Then, we employ the well-known Runge-Kutta method to solve this system, which aids in addressing the nonlinearity and the lack of data on $\partial \Omega \setminus \Gamma$.
This new approach is called the dimensional reduction method.
By converting the high-dimensional problem into a 1D problem, we achieve exceptional computational speed. Numerical results are provided to support the effectiveness of our approach.
\end{abstract}



\section{Introduction}

Let $d \geq 1$ be the spatial dimension.
Let $a < b$, $a_2 < b_2$, $\dots,$ $a_d < b_d$ be constants.
Define 
\[
	\widetilde \Omega =
	\left\{
		\begin{array}{ll}
	 		\ds\prod_{i=2}^d (a_i,b_i)  & d > 1,\\
			\emptyset &d = 1,
		 \end{array}
	\right.
	\quad
	\mbox{and}
	\quad
	\Omega = (a, b) \times \widetilde \Omega.
\]
For points $\x \in \Omega,$ we write $\x = (x, \widetilde \x)$ where $x \in (a, b)$ is the first entry of $\x$ and $\widetilde \x = (x_2, \dots, x_d) \in \widetilde \Omega$  consists of the $2^{\rm nd}$, $\dots$, $d^{th}$ entries of $\x.$
 Let $T$ be a positive number representing the final time. 
 Throughout the paper, we denote by $\Omega_T$ and $\widetilde \Omega_T$ the sets $\Omega \times (0, T)$ and $\widetilde \Omega \times (0, T)$ respectively.
 Let 
\begin{equation}
	\Gamma = \big\{(a, \widetilde \x): \widetilde \x \in \widetilde \Omega\big\}
	\quad
	\mbox{and}
	\quad
	\Gamma_T = \Gamma \times (0, T)
\end{equation}
represent the measurement site.
Note that $\Gamma$ is the left side of $\partial \Omega$ which is perpendicular to the $x-$axis.
Let $u$ be a $C^2$ solution to the following nonlinear parabolic equation
\begin{equation}
	u_t(\x, t) = \Delta u(\x, t) + F\left(\x, t, u(\x,t),\nabla u(\x,t)\right),
	\quad (\x,t) \in \Omega_T
	\label{main_eqn}
\end{equation}
where the nonlinearity
$
	F: \Omega \times \mathbb{R} \times \mathbb{R} \times \mathbb{R}^d \to \mathbb{R}
$ 
is a smooth function.
We are interested in the following nonlinear version of the inverse heat conduction problem:
\begin{Problem}[Inverse Heat Conduction Problem]
Given the  measurement on $\Gamma_T$
\begin{equation}
g(\widetilde \x,t) = u(a, \widetilde \x,t) \mbox{ and } q(\widetilde \x,t)  = \frac{\partial}{\partial x} u(a, \widetilde \x,t) \quad \mbox{ for } (\widetilde \x,t) \in \widetilde \Omega_T,
\label{1.3}
\end{equation}
determine the function $u(\x,t)$ for $(\x,t) \in \Omega_T$.
\label{IHCP}
\end{Problem}

Assume that the nonlinearity $F$ satisfies the following growth condition
\begin{equation*}
	F(\x, t, s, \p) \leq C\big(
		|s| + |\p| + |f|
	\big)
\end{equation*}
for all $\x \in \overline \Omega$, $t \in [0, \infty),$ $s \in \R$, and $\p \in \R^d$ for some function $f \in L^2(\Omega),$ and for some positive constant $C$.
Then, the uniqueness of Problem \ref{IHCP} is a direct consequence of Theorem 2.2.1 and Theorem 2.3.1 in \cite{KlibanovLiBook}. 
In practice, the solution $u(\x, t),$ $(\x, t) \in \Omega_T$ for Problem \ref{IHCP} could represent temperature or pollutant distribution. As a result, computing $u$ has significant value in numerous real-world applications, such as determining the spatial temperature distribution within an inaccessible domain $\Omega$, given external measurements of heat and heat flux, \cite{Klibanov:ip2006}; identifying pollutant levels on the river or lake surfaces \cite{BadiaDuong:jiip2002}; and effectively monitoring heat conduction processes in steel, glass, and polymer-forming industries, as well as nuclear power plans \cite{LiYamamotoZou:cpaa2009}. 
Due to its significance, the inverse heat conduction problem was studied intensively.
We cite \cite{Nehad1998, alifanov-book1,alifanov-book2, Arora:2019,beck-book,Burggraf1964,Ginsberg1963, Felde2014, isakov, Hao1992,Hao1994,haob, JARNY19912911, Loulou:2006, Mohebbi2021,Ozisik2017, Roy2023} and references therein. for other studies including theoretical and numerical studies about the related versions of the inverse heat conduction problem.
Moreover, since we do not impose a special structure on the nonlinearity $F$,  
equation \eqref{main_eqn} is not limited to be the model for heat transferring.
As a result, Problem \ref{IHCP} has applications in some other fields.
For example, when $F$ takes the form $u(1 - u)$ (or more generally, $F(u) = u(1 - |u|^{\alpha})$ for some $\alpha > 0$), the parabolic equation \eqref{main_eqn} becomes the higher-dimensional version of the well-known Fisher (or Fisher-Kolmogorov) equation \cite{Fisher:ae1937}.
 It is important to mention that the Fisher equation naturally appears in mathematical models found in ecology, physiology, combustion, crystallization, plasma physics, and general phase transition problems, as seen in \cite{Fisher:ae1937}. Hence, finding solutions to nonlinear parabolic equations is of considerable practical significance. This problem is worth investigating.

As the primary objective of this paper is to determine a solution for \eqref{main_eqn} from the Cauchy data observed on $\Gamma_T$ \eqref{1.3}, we assume that the solution exists. For the sake of completeness, we recall here a collection of conditions that ensure the existence. In the theory of PDEs, in order to determine a solution to a parabolic equation like \eqref{main_eqn}, one needs to impose an initial condition and a boundary condition. Assume that the initial condition 
 $u(\mathbf{x}, 0)$ belongs to $H^{2+\beta}(\mathbb{R}^d)$ for some $\beta \in [0, 1 + 4/d]$, and for all $\mathbf{x} \in \Omega$, $t \in [0, T]$, $s \in \mathbb{R}$, and $\mathbf{r} \in \mathbb{R}^d$, we have 
\[\left|F(\mathbf{x}, t, s, \mathbf{p})\right| \leq C\max\{(1 + |\mathbf{p}|)^2, 1 + |s|\}\]
for some positive constant $C$. 
Assume further that the boundary value of the function $u$ is sufficiently smooth.
Then, due to Theorem 6.1 in \cite[Chapter 5, \S 6]{LadyZhenskaya:ams1968} and Theorem 2.1 in \cite[Chapter 5, \S 2]{LadyZhenskaya:ams1968}, equation \eqref{main_eqn}, together with suitable initial condition, has a unique solution with $|u(\mathbf{x},t)|\leq M_1$, $|\nabla u(\mathbf{x},t)| \leq M_2$ for some positive constants $M_1,M_2$, and $u(\mathbf{x},t) \in H^{2+\beta, 1+\beta/2}(\mathbb{R}^d \times [0,T])$.

Numerically solving Problem \ref{IHCP} is highly challenging. First, the problem is severely ill-posed in the sense that a minor amount of noise can lead to substantial computational errors, see a brief discussion in Section \ref{sec 3}. 
Another difficulty lies in the nonlinear behavior of Problem \ref{IHCP}. To solve nonlinear problems numerically, one employs methods based on least-squares optimization including a Tikhonov regularization term. The cost functional can be defined in numerous ways. One common example of a cost functional is:
\begin{equation}
	J(u) = \|u_t(\x, t) - \Delta u(\x, t) - F(\x, t, u(\x, t), \nabla u(\x, t))\|^2_{L^2(\Omega_T)} + \epsilon \|u(\x, t)\|_H^2.
	\label{1.6}
\end{equation}
Here, the regularization parameter $\epsilon$ and the regularization norm $\|\cdot\|_H$ are chosen by the users.
To obtain a numerical solution, one minimizes $J$ subject to the constraints $u(a, \widetilde \x, t) = g(\widetilde \x, t)$ and $\partial_{x}u(a, \x, t) = q(\widetilde \x, t)$ for $(\widetilde \x, t) \in \widetilde \Omega_T.$
The nonlinearity $F$ may cause the functional $J$ to be nonconvex and to have multiple local minima. Therefore, ensuring the successful identification of the global minimizer of $J$ requires a good initial guess. However, in practice, this requirement is often too stringent, as an initial guess may not always be available in real-world applications.
Regarding numerical methods that do not request good initial guesses of the true solutions, we draw the reader's attention to several recent works \cite{AbhishekLeNguyenKhan, LeCON2023, LeNguyen:jiip2022} in which problems similar to Problem \ref{IHCP} were addressed.
Although these papers focus on determining the initial stages for general nonlinear parabolic equations,
which differ from that of Problem \ref{IHCP}, numerical methods in these papers allow for the computation of all stages of the solution. It is worth noting that, unlike Problem \ref{IHCP}, the approaches in \cite{AbhishekLeNguyenKhan, LeCON2023, LeNguyen:jiip2022} require full boundary measurements, whereas, in the present paper, we only need partial data to determine the solution. 
We would like to highlight several recent studies \cite{KlibanovNguyenTran:JCP2022, LeNguyen:JSC2022, LeNguyenTran:CAMWA2022, NguyenKlibanov:ip2022} for numerical methods based on Carleman estimates to compute the solutions of nonlinear elliptic equations, and the paper \cite{NguyenKlibanov:ip2022} that focuses on solving nonlinear hyperbolic equations. Like the works \cite{AbhishekLeNguyenKhan, LeCON2023, LeNguyen:jiip2022}, the data required in \cite{KlibanovNguyenTran:JCP2022, LeNguyen:JSC2022, LeNguyenTran:CAMWA2022, NguyenKlibanov:ip2022} are obtained from measurements on the entire boundary of the computational domain.

In this paper, we introduce a novel approach for addressing Problem \ref{IHCP}, called the dimensional reduction method.
The rationale behind this name lies in the center of our approach, which transforms Problem \ref{IHCP} into the problem of solving a system of ordinary differential equations (ODEs) defined along the $x-$axis.
Let $\{\Psi_n\}_{n \geq 1}$ be an orthonormal basis of $L^2(\widetilde \Omega \times (0, T))$. For each $\x = (x, \widetilde \x) \in \Omega$ with $a \leq x \leq b$ and $\widetilde \x \in \widetilde \Omega,$ we approximate $u(x, \widetilde \x, t)$ by
\begin{equation}
	u(x, \widetilde \x, t) 
	= \sum_{n = 1}^{\infty} u_n(x) \Phi_n(\widetilde \x, t) 
	\approx u(x, \widetilde \x, t) = \sum_{n = 1}^{N} u_n(x) \Phi_n(\widetilde \x, t) 
	\label{1.6666}
\end{equation}
for some cut-off number $N$.
Here, 
\[
	u_n(x) = \int_{\widetilde \Omega_T} u(x, \widetilde \x, t) \Phi_n(\widetilde \x,t) d\widetilde \x dt,
	\quad n \geq 1.
\]
We then substitute the approximation \eqref{1.6666} into \eqref{main_eqn} to derive a second-order ODE system for $\{u_n(x)\}_{n = 1}^N$ with Cauchy data provided at the endpoint $a$. Using a standard variable change, as outlined in any ODE textbook, we can convert this system into a set of first-order equations. Classical Euler or Runge-Kutta methods can then solve the new system.
We have reduced the $(d + 1)-$dimensional problem to a one-dimensional problem.
We name our approach the {\it dimensional reduction} method.
It possesses two significant strengths
\begin{enumerate}
	\item It rapidly delivers solutions;
	\item  It does not request a good initial guess of the true solution.
\end{enumerate}
On the other hand, the dimensional reduction approach has a limitation. 
It can only compute the solution $u(\x, t)$ when $\x$ is in a domain near the measurement site $\Gamma$.
 The performance declines as $\x$ gets further away from $\Gamma$. However, this drawback is acceptable because data is not completely provided on the whole $\partial \Omega$. 
A crucial factor in the success of the dimensional reduction method is the selection of an appropriate basis $\{\Phi_n\}_{n \geq 1}$. In this paper, we employ the polynomial-exponential basis, introduced in \cite{Klibanov:jiip2017}
and further developed in \cite{NguyenLeNguyenKlibanov:arxiv2023}. It is important to note that employing well-known bases like Legendre polynomials or trigonometric functions in the widely-used Fourier expansion may not be optimal. For a comparison of the performance of these bases in computation, we direct readers to \cite{LeNguyenNguyenPowell:JOSC2021, NguyenLeNguyenKlibanov:arxiv2023}.

The paper is organized as follows. 
In Section \ref{sec 2}, we present our dimensional reduction method to solve Problem \ref{IHCP}.
In Section \eqref{sec 3}, we outline some heuristic discussion about the efficiency of the dimensional reduction method. 
In Section \ref{sec numerical study}, we present the implementation.
Section \ref{sec 5} is for numerical results.
Section \ref{sec remarks} is for some concluding remarks.

\section{The dimensional reduction approach} \label{sec 2}

We first recall the polynomial-exponential basis that plays a crucial role to our approach.
Let $a < b$ be two real numbers. 
For each positive integer number $n$, define $\phi_n(s) = s^{n-1} e^{s - (a+b)/2}$ for all $s \in (a, b).$
The set $\{\phi\}_{n \geq 1}$ is complete in $L^2(a, b)$. 
Applying the Gram-Schmidt orthonormalization process on this set, we obtain an orthonormal basis, denoted by $\big\{\Psi_n^{(a, b)}\big\}_{n \geq 1}$ of $L^2(a, b)$. 
This basis was introduced by Klibanov in \cite{Klibanov:jiip2017} to solve inverse problems.
It plays a crucial role in the center of our reduction dimension method.
Recall that the computational domain $\Omega$ is the cube $(a, b) \times \ds\prod_{i = 2}^d (a_i, b_i)$, $\widetilde \Omega = \ds\prod_{i = 2}^d (a_i, b_i)$, and $T > 0$ is the final time.
For any multi-index $\bn = (n_2, \dots, n_d, n_t) \in \mathbb{N}^d$, we define 
\[
	P_{\bn}(x_2, \dots, x_d, t) = \Psi_{n_t}^{(0, T)}(t)\prod_{i = 2}^d \Psi_{n_i}^{(a_i, b_i)}(x_i).
\]
It is obvious that $\{P_{\bn}\}_{\bn \in \R^d}$ is an orthonormal basis of $L^2(\widetilde \Omega \times (0, T)).$
We name it the {\it polynomial-exponential basis}.
The key step of our method is to truncate the Fourier expansion with respect to the basis $\{P_\bn\}_{\bn \in \R^d}$.
The procedure is described below.

The Fourier expansion of $u$ at $(\x, t) = (x, \widetilde \x, t) \in \Omega_T$ with $(\widetilde \x, t) \in \widetilde \Omega_T$
is given by
\begin{equation}
	u(\x, t) = u(x, \widetilde \x, t) = \sum_{\bn \in \R^d} u_\bn(x) P_{\bn}(\widetilde \x, t) 
\label{2.1}
\end{equation}
where
\begin{equation}
	u_{\bn}(x) = \int_{\widetilde \Omega_T} u(x, \widetilde \x, t) P_{\bn}(\widetilde \x, t) d\widetilde \x dt.
	\label{2.2}
\end{equation}
To transform the high-dimensional Problem \ref{IHCP} into a 1D problem along the $x-$axis, we approximate the Fourier expansion \eqref{2.1} by truncating it as follows:
\begin{equation}
	u(\x, t) = u(x, \widetilde \x, t) 
	\approx \sum_{\bn \in {\bf N}} u_{\bn}(x) P_{\bn}(\widetilde \x, t),
	\label{2.3333}
\end{equation}
where
\[
	\bN = \{(n_2, \dots, n_d, n_t): 1 \leq n_2 \leq N_2, \dots ,1 \leq n_d \leq N_d, 1 \leq n_t \leq N_t\}
\] and
the cut-off numbers $N_2, \dots, N_d,$ and $N_t$ will be determined later based on the provided data. We will discuss this ``data-driven" approach in Section \ref{sec numerical study}.

Plugging approximation \eqref{2.3333} into \eqref{main_eqn} gives
\begin{multline}
 \sum_{\bn \in {\bf N}} u_{\bn}(x) \frac{\partial}{\partial t}P_{\bn}(\widetilde \x, t) 
 = 
  \sum_{\bn \in {\bf N}} u_{\bn}''(x) P_{\bn}(\widetilde \x, t)
+  \sum_{\bn \in {\bf N}} u_{\bn}(x) \Delta_{\widetilde \x}P_{\bn}(\widetilde \x, t)
\\
 + F\left(\x, t,  \sum_{\bn \in {\bf N}} u_{\bn}(x) P_{\bn}(\widetilde \x, t), \sum_{\bn \in {\bf N}} u_{\bn}'(x) P_{\bn}(\widetilde \x, t),  \sum_{\bn \in {\bf N}} u_{\bn}(x) \nabla_{\widetilde \x} P_{\bn}(\widetilde \x, t)\right)
 \label{2.3}
\end{multline}
for all $(x, \widetilde \x, t) \in \Omega_T.$
Here,
\begin{equation*}
	\Delta_{\widetilde \x}P_{\bn}(\widetilde \x, t) = 
	\sum_{i = 2}^d \frac{\partial^2}{\partial x_i^2} P_{\bn}(\widetilde \x, t),
	\quad
	\mbox{and}
	\quad
	\nabla_{\widetilde \x}P_{\bn}(\widetilde \x, t) = \left(
			\frac{\partial}{\partial x_i} P_{\bn}(\widetilde \x, t)
		\right)_{i = 2}^d 
\end{equation*}
\begin{Remark}
	The derivation of equation \eqref{2.3} may not be entirely theoretically rigorous due to the term-by-term differentiation of the Fourier expansion in \eqref{2.1}. The system of ODEs in \eqref{2.3} should be viewed as an approximation context of Problem \ref{IHCP}. Analyzing the behavior of \eqref{2.3} when $\bN \to \mathbb{N}^d$, meaning all cut-off numbers $N_2, \dots, N_d, N_t$ tending to $\infty$, is highly challenging. This issue is not addressed in this paper, which focuses on computation.
However, term-by-term differentiation of the Fourier expansion using the polynomial-exponential basis has been numerically validated in \cite{NguyenLeNguyenKlibanov:arxiv2023}. It was demonstrated in \cite{NguyenLeNguyenKlibanov:arxiv2023} that calculating up to second-order derivatives by term-by-term differentiating the truncated Fourier expansion with respect to the polynomial-exponential basis yields accurate results even with highly noisy data. This finding somewhat supports the accuracy of \eqref{2.3}.
\end{Remark}
For every $\bm \in \bN$, multiplying $P_\bm(\widetilde \x, t)$ into both sides of \eqref{2.3} and then integrating the obtained equation on $\widetilde \Omega_T$ gives
\begin{multline}
\int_{\widetilde \Omega_T} \sum_{\bn \in {\bf N}} u_{\bn}(x) \frac{\partial}{\partial t}P_{\bn}(\widetilde \x, t) P_\bm(\widetilde \x, t) d\widetilde \x d t
\\
 = 
  \int_{\widetilde \Omega_T}\sum_{\bn \in {\bf N}} u_{\bn}''(x) P_{\bn}(\widetilde \x, t) P_\bm(\widetilde \x, t) d\widetilde \x d t
+  \int_{\widetilde \Omega_T}\sum_{\bn \in {\bf N}} u_{\bn}(x) \Delta_{\widetilde \x}P_{\bn}(\widetilde \x, t) P_\bm(\widetilde \x, t) d\widetilde \x d t
\\
 + \int_{\widetilde \Omega_T} F\left(\x, t,  \sum_{\bn \in {\bf N}} u_{\bn}(x) P_{\bn}(\widetilde \x, t), \sum_{\bn \in {\bf N}} u_{\bn}'(x) P_{\bn}(\widetilde \x, t),  \sum_{\bn \in {\bf N}} u_{\bn}(x) \nabla_{\widetilde \x} P_{\bn}(\widetilde \x, t)\right) P_\bm(\widetilde \x, t) d\widetilde \x d t.
 \label{2.4}
\end{multline}
Equation \eqref{2.4} is reduced to
\begin{equation}
	u_\bm''(x) = \sum_{\bn \in {\bf N}} s_{\bm \bn} u_{\bn}(x) - 
	F_{\bm}\left(x, \{u_{\bn}(x)\}_{\bn \in \bN}, \{u_{\bn}'(x)\}_{\bn \in \bN}\right),
	\label{2.6}
\end{equation}
for $x \in (a, b)$,
where
\begin{equation}
	s_{\bm \bn} = \sum_{\bn \in {\bf N}} \int_{\widetilde \Omega_T}   \Big[\frac{\partial}{\partial t}P_{\bn}(\widetilde \x, t)
	- \Delta_{\widetilde \x}P_{\bn}(\widetilde \x, t)
	\Big] P_\bm(\widetilde \x, t) d\widetilde \x d t,
\end{equation}
and
\begin{multline}
	F_{\bm}\left(x, \{u_{\bn}(x)\}_{\bn \in \bN}, \{u_{\bn}'(x)\}_{\bn \in \bN}\right) 
	\\
	= \int_{\widetilde \Omega_T} F\left(\x, t,  \sum_{\bn \in {\bf N}} u_{\bn}(x) P_{\bn}(\widetilde \x, t), \sum_{\bn \in {\bf N}} u_{\bn}'(x) P_{\bn}(\widetilde \x, t),  \sum_{\bn \in {\bf N}} u_{\bn}(x) \nabla_{\widetilde \x} P_{\bn}(\widetilde \x, t)\right) P_\bm(\widetilde \x, t) d\widetilde \x d t.
\end{multline}
The Cauchy data for $\{u_{\bm}\}_{\bm \in \bN}$ can be computed from the provided data of Problem \ref{IHCP} and formula \eqref{2.2}
\begin{align}
	u_{\bm}(a) &= g_\bm :=  \int_{\widetilde \Omega_T} g(a, \widetilde \x, t) P_{\bn}(\widetilde \x, t) d\widetilde \x dt,
	\label{2.9}
	\\
	u_{\bm}'(a) &= q_\bm := \int_{\widetilde \Omega_T} q(a, \widetilde \x, t) P_{\bn}(\widetilde \x, t) d\widetilde \x dt.
	\label{2.10}
\end{align}

By \eqref{2.6}, \eqref{2.9}, and \eqref{2.10}, we have reduced Problem \ref{IHCP} to the problem of computing $\{u_{\bm}(x)\}_{\bm \in \bN}$ that satisfies the initial value problem
\begin{equation}
	\left\{
		\begin{array}{rcll}
			u_\bm''(x) &=&
		\ds \sum_{\bn \in {\bf N}} s_{\bm \bn} u_{\bn}(x) - 
	F_{\bm}\left(x, \{u_{\bn}(x)\}_{\bn \in \bN}, \{u_{\bn}'(x)\}_{\bn \in \bN}\right) &x \in (a, b)\\
	u_{\bm}(a) &=& g_{\bm}
	\\
	u'_\bm(a) &=&q_\bm 
		\end{array}
	\right.
	\label{2.11}
\end{equation}

We opt for the Runge-Kutta method to solve \eqref{2.11} since it is more stable compared to the classical Euler method. To apply the Runge-Kutta method to \eqref{2.11}, we must rewrite \eqref{2.11} to a system of first-order equations. To do this, we define
\begin{equation}
	w_{\bN}(x) = \left[
		\begin{array}{c}
			\{u_{\bm}(x)\}_{\bm \in \bN}
			\\
			\{u_{\bm}'(x)\}_{\bm \in \bN}
		\end{array}
	\right]
	\quad
	\mbox{for all } x \in (a, b).
	\label{2.1212}
\end{equation}
It follows from \eqref{2.11} that
\begin{align}
	w_{\bN}'(x) &= \left[
		\begin{array}{c}
			\{u_{\bm}'(x)\}_{\bm \in \bN}
			\\
			\{u_{\bm}''(x)\}_{\bm \in \bN}
		\end{array}
	\right]
	\nonumber
	\\
	&= 
	\left[
		\begin{array}{l}
			\{u_{\bm}'(x)\}_{\bm \in \bN}
			\\
			\Big\{\ds \sum_{\bn \in {\bf N}} s_{\bm \bn} u_{\bn}(x) - 
	F_{\bm}\left(x, \{u_{\bn}(x)\}_{\bn \in \bN}, \{u_{\bn}'(x)\}_{\bn \in \bN}\right)\Big\}_{\bm \in \bN}
		\end{array}
	\right]	\label{2.12}
\end{align}
for $x \in (a, b).$
Define ${\bf F}(x, w_{\bN}(x))$ as the expression in \eqref{2.12}.
Then, by using \eqref{2.11}, \eqref{2.1212}, and \eqref{2.12}, we obtain 
\begin{equation}
	\left\{
		\begin{array}{rcll}
			w_{\bN}'(x) &=& {\bf F}(x, w_{\bN}(x)) &x \in (a, b)\\
			w_{\bN}(a) &=& 
			\left[
				\begin{array}{c}
					\{g_\bm\}_{\bm \in \bN}\\
					\{q_\bm\}_{\bm \in \bN}
				\end{array}
			\right].
		\end{array}
	\right.
	\label{2.14}
\end{equation}

\begin{Remark}
The center of the dimension reduction approach lies in the derivation of \eqref{2.14}, which represents a system of first-order ODEs along the $x-$axis.
 Although \eqref{2.14} is nonlinear, solving it is relatively straightforward. Various methods can be employed, such as the Euler and Runge-Kutta methods, which are covered in ODE textbooks. Once the solution $\{w_\bm(x)\}_{\bm \in \bN}$, and hence $\{u_\bm(x)\}_{\bm \in \bN}$, $x \in (a, b)$ is obtained, the function $u(\x, t)$ for $(\x, t) \in \Omega_T$ can be computed using \eqref{2.3333}.
 
The primary concern with this approach is that solving the initial value problem \eqref{2.14} may be unstable. However, this is understandable given that the problem being examined is severely ill-posed, as demonstrated by an explicit example in Section \ref{sec 3}. 
This example indicates that computing the solution to Problem \ref{IHCP} might not be feasible when $x > a$.
 Nevertheless, our dimension reduction approach can provide solutions within an open layer near the measurement site $\Gamma$. We will further discuss the reasons for the success of this approach in Section \ref{sec 3} and Section \ref{sec numerical study}.
\label{rem 2.2}
\end{Remark}

\section{Some discussions about the dimensional reduction method}\label{sec 3}

We will discuss the difficulties associated with numerically solving Problem \ref{IHCP} and explain why the dimension reduction method is an appropriate choice for addressing these challenges.

To demonstrate the ill-posedness of Problem \ref{IHCP},
we examine the case $d = 1$ and $F = 0$. 
Let $u^*$ be a solution to the following equation 
\begin{equation}
	u_t(x, t) = u_{xx}(x, t),
	\quad (x, t) \in (a_1, b_1) \times (0, T)
	\label{1.4}
\end{equation}
Let $g^*(t)$ and $q^*(t)$ represent the exact measurements of the function $u^*(0, t)$ and $u^*_x(0, t)$ respectively. 
Assume that the measurements are  given by
\[
	g(t) = g^*(t) + \eta_1(t),
	\quad
	q(t) = q^*(t) + \eta_2(t)
\]
where the small noises $\eta_1$ and $\eta_2$ are given by
\begin{equation}
	\eta_1(t) = 2n^{-2} 
		\sin(2n^2 t),
	\quad
	\eta_2(t) = 0.
	\label{2.2}
\end{equation}
Here, $n$ is a positive integer.
Then, one can directly verify that the function  
\[
	u(x, t) = u^*(x, t) + n^{-2} [
		e^{nx} \sin(2n^2 t + nx) + e^{-n x} \sin(2n^2 t - nx)
	],
	\quad
	x \in \R, t \in (0, \infty)
\] solves Problem \eqref{IHCP}.
The error in computation is 
\[
	u(x, t) - u^*(x, t) = n^{-2} [
		e^{nx} \sin(2n^2 t + nx) + e^{-n x} \sin(2n^2 t - nx)
	],
\]
which becomes large when when $x > 0$ and $n \gg 1.$
In other words, it is not possible to approximate the true solution $u^*(x, t)$ unless a special technique is employed.
This special technique in this paper is the truncation of all highly oscillated components of the data, and hence the noise.
This technique is demonstrated in the approximation \eqref{2.3333}.
More precisely,
the noise $\eta_1$ in example \eqref{2.2} above emphasizes that the large error could be due to the high oscillation.
 As a result, in order to reduce computational errors, we propose the elimination of highly oscillatory components from the given data. To accomplish this, we truncate all of the corresponding highly oscillatory terms in the Fourier expansion of the solution $u$, as described in equation \eqref{2.3333}. By removing these highly oscillatory terms, the resulting 1D approximation model in equation \eqref{2.14} exhibits a ``more regular" behavior compared to the original model.
The numerical evidence for this observation will be presented in Section \ref{sec 5} in which we can compute the solution in a layer with positive volume near the measurement site $\Gamma.$ 
We would like to note that since the initial value problem in \eqref{2.3333} is not ``globally" stable, the reduction dimensional method does not completely remove the ill-posedness.
When $\x$ is too far from $\Gamma$, our method does not provide a reliable solution.

On the other hand, our truncation technique, or the removal of highly oscillatory terms, is somewhat consistent with the Tikhonov regularization, which is a conventional method to address the ill-posedness.
We consider that Tikhonov regularization is a great approach to approximate an ill-posed problem with another well-posed problem. 
This can be done by adding a Tikhonov regularization term into a cost functional, as shown in equation \eqref{1.6}, when the functional space $H$ possesses a strong norm. Possible choices for $H$ include $H^2(\Omega_T),$ $H^3(\Omega_T),$ or even $H^p(\Omega_T),$ where $p$ is a large integer. 
Heuristically, including these norms in the cost functional prevents the explosive behavior of the computed solution's derivatives, leading to reduced oscillation in the computed solution.
Since our truncation technique in \eqref{2.3333} has the same effect of reducing the high oscillation, we conclude that when we cut-off the series in \eqref{2.1}, we somewhat regularize the original ill-posed problem.

We next present the significance of the polynomial-exponential basis used in this paper. One may wonder why we selected this basis, among an infinite number of others, for the Fourier expansion in \eqref{2.1}. The reason is that popular bases, such as Legendre polynomials or trigonometric functions, may not be appropriate. This is because the initial function $P_{\bm_0} = P_{(1, \dots, 1)} $ of these bases is a constant, resulting in an identically zero derivative. Consequently, the corresponding Fourier coefficient $u_{\bm_0}(x)$ in the sums $ \sum_{\bn \in {\bf N}} u_{\bn}(x) \frac{\partial}{\partial t}P_{\bn}(\widetilde \x, t)$ and $ \sum_{\bn \in {\bf N}} u_{\bn}(x) \Delta_{\widetilde \x}P_{\bn}(\widetilde \x, t)$ in \eqref{2.3} is overlooked, leading to reduced accuracy. The polynomial-exponential basis used in this paper is appropriate for \eqref{2.3} since it satisfies the necessary condition that derivatives of $P_{\bn}$, $\bn \in \mathbb{N}^d$, are not identically zero. We have demonstrated the effectiveness of the polynomial-exponential basis in various studies, including \cite{LeNguyenNguyenPowell:JOSC2021, NguyenLeNguyenKlibanov:arxiv2023}. In \cite{LeNguyenNguyenPowell:JOSC2021}, we used the polynomial-exponential and traditional trigonometric bases to expand wave fields and solve problems in photo-acoustic and thermo-acoustic tomography. Results indicated that the polynomial-exponential basis performed much better. In \cite{NguyenLeNguyenKlibanov:arxiv2023}, we employed Fourier expansion to compute derivatives of data interrupted by noise. Results showed that the polynomial-exponential basis was more accurate than the trigonometric basis when differentiating term-by-term of the Fourier expansion in solving ill-posed problems. Therefore, the polynomial-exponential basis should be used when term-by-term differentiation of the Fourier expansion is required.

\section{Numerical implementation} 
\label{sec numerical study} 

Recall that Section \ref{sec 2} outlines the derivation of the 1D approximation initial value problem \eqref{2.14}, which forms the center of the dimensional reduction method. As discussed in the first paragraph of Remark \ref{rem 2.2}, computing the solution to \eqref{2.14} enables us to directly obtain the solution to Problem \ref{IHCP}. For ease of use, we summarize the steps to solve Problem \ref{IHCP} in Algorithm \ref{alg}.

\begin{algorithm}[!ht]
\caption{\label{alg}The dimensional reduction method for the nonlinear inverse heat conduction problem}
	\begin{algorithmic}[1]
	\State \label{s1} Choose the cut-off numbers $N_1, N_2, \dots, N_t$. Define 
	\begin{equation}
		\bN = \{(n_1, \dots, n_d, n_t): 1 \leq n_i \leq N_i,  1 \leq i \leq d,1 \,\mbox{and}\, \leq n_t \leq N_t\}.
		\label{4.1}
	\end{equation}
	\State \label{s2} Compute the initial data $\{g_\bm\}_{\bm \in \bN}$ and $\{q_\bm\}_{\bm \in \bN}$ by formulas \eqref{2.9} and \eqref{2.10} respectively.
	\State \label{s3} Solve the system of ODEs \eqref{2.14} with initial values computed in the previous step.
	Denote by $w_\bN^{\rm comp}$ the computed solution.
	\State \label{s4} Extract the knowledge of $\{u_{\bm}^{\rm comp}\}_{\bm \in \bN}$ using the definition of  $w_\bN^{\rm comp}$ in \eqref{2.1212}.
	\State \label{s5} Compute the desired solution $u^{\rm comp}(\x,t), \x \in \Omega,t \in [0,T]$ using \eqref{2.3333} with $u_{\bn}$ being replaced by $u^{\rm comp}_{\bn}.$
\end{algorithmic}
\end{algorithm}

For simplicity, we implement and present numerical results due to Algorithm \ref{alg} in one and two dimensions, i.e., $d=1$ and $d=2$. 
In these tests, we will choose the parabolic equations of the form \eqref{main_eqn} that we already know the analytic solutions.
The noiseless simulated data $g^*$ and $q^*$ are extracted by those solutions. 
For some noise level $\delta > 0$, we generate the noisy data $g^\delta$ and $q^\delta$ by the formulas
\begin{equation}
	g^{\delta} = g^*(1 + \delta {\rm rand})
	\quad
	\mbox{and}
	\quad
	q^{\delta} = q^*(1 + \delta {\rm rand}),
\end{equation}
where rand is a function taking uniformly distributed random numbers in $[-1, 1].$
	
\begin{Remark}
	We observe that Algorithm \ref{alg} can produce reliable solutions to Problem \ref{IHCP} over the entire domain $\Omega \times (0,T)$ only when the data is noise-free. However, when Algorithm \ref{alg} is tested with noisy data, we find that for noise levels of $\delta = 5\%$ and $\delta = 10\%$, the algorithm yields good numerical solutions $u^{\rm comp}(\x,t)$ for $\x \in \Omega$ and $t < T'$, where $T'$ is a number less than $T$. This issue may stem from the instability of solving initial value systems of ODEs. For our numerical study, we set $T=1.5$ and $T'=1.3$. Accordingly, we only present the numerical results over the domain $\Omega \times (0,T')$. This approach is related to the so-called "future time method" by Beck \cite{beck-book}: use more Cauchy data to stabilize the solution to the inverse heat conduction problem.  
\end{Remark}


We use the finite difference scheme to implement Algorithm \ref{alg}. We divide the interval $[0,T]$ into a uniform partition as
\[
	\mathcal T = \{t_i = (i-1) d_t: i = 1, \dots, N_T\}
\] where $N_T$ is a positive integer and $d_t = T/(N_T - 1).$ 
For all computational results presented in this paper, we set $T = 1.5$ and $N_T = 601$.
We also discretize the domain $\Omega$ with a uniform mesh. 
\begin{enumerate}
	\item For $d= 1$, we set $a = 0$ and $b = 1$; i.e., $\Omega = (0, 1)$ and compute the solution on 
	\[
		\mathcal{G}_{\rm 1d} = \{x_i = (i - 1)d_x: i = 1, \dots, N_x\}
	\]
	where $N_x$ is a positive integer and $d_x = (b - a)/(N_x - 1).$ In this paper, $N_x = 201$.
	\item For the case where $d = 2$, we assign the values $a = 0$, $b = 0.2$, $a_2 = -1$, and $b_2 = 1$ to define the domain $\Omega = (0, 0.5) \times (-1, 1)$. We calculate the solution on the grid $\mathcal{G}_{\rm 2d}$ defined as
	\[
		\mathcal{G}_{\rm 2d} = \big\{
			(x_i = (i - 1)d_{x}, y_j = -1 + (j - 1)d_{y}): i = 1, \dots, N_{x}, j = 1, \dots, N_y, 
		\big\}
	\] 
	where $N_x$ and $N_y$ are positive integers and the step sizes in the $x-$direction and the $y-$direction are $d_x = (b - a)/(N_x-1)$ and $d_y = (b_2 - a_2)/(N_y - 1)$ respectively. For the present study, we select $N_x = 201$ and $N_y = 1001$.
\end{enumerate}


The implementation of steps \ref{s2}, \ref{s4}, and \ref{s5} of Algorithm \ref{alg} are straightforward. 
We only present details for Steps \ref{s1} and \ref{s3}.

\subsection{ Step \ref{s1} of Algorithm \ref{alg}}
We use a data-driven approach for choosing cut-off numbers.
More precisely, the knowledge of given data $g(\x, t)$, $(\x, t) \in \Gamma_T$ is helpful in determining the cut-off numbers for step \ref{s1} of Algorithm \ref{alg}. The procedure is based on a trial-and-error process. 
Inspired by \eqref{2.3333} and the fact that $u(a_1, \widetilde \x, t) = g(a_1, \widetilde \x, t)$ for all $(\x', t) \in \widetilde \Omega_T$, for each $N_2, \dots, N_d, N_t \geq 1$, we define the function $\varphi: \mathbb{N}^d \to \R$  that represents the relative distance between the data $g$ and the truncation of its Fourier expansion with respect to the polynomial-exponential basis.
The function $\varphi$ is defined as follows
\begin{equation}
	\varphi(\bN)= \frac{\|g- \sum_{\bn \in \bN} g_{\bn} P_{\bn}(\widetilde \x, t)\|_{L^{\infty}(\Gamma_T)}}{\|g\|_{L^{\infty}(\Gamma_T)}},
	\quad
	\mbox{where}\,
	g_{\bn} = \int_{\Omega'_T} g(a_1, \widetilde \x, t )P_{\bn}(\widetilde \x, t)d\widetilde \x dt
	\label{4.2}
\end{equation}
where $\bN = \{ (n_1, \dots, n_d, n_t): 1 \leq n_i \leq N_i,  1 \leq i \leq d,1 \,\mbox{and}\, \leq n_t \leq N_t\}$ as in \eqref{4.1}.
 We increase the numbers $N_2, \dots, N_d, N_t$ until $\eta_\bN$ is sufficiently small.
 
\begin{figure}[h!]
\centering
	\subfloat[]{\includegraphics[width=.3\textwidth]{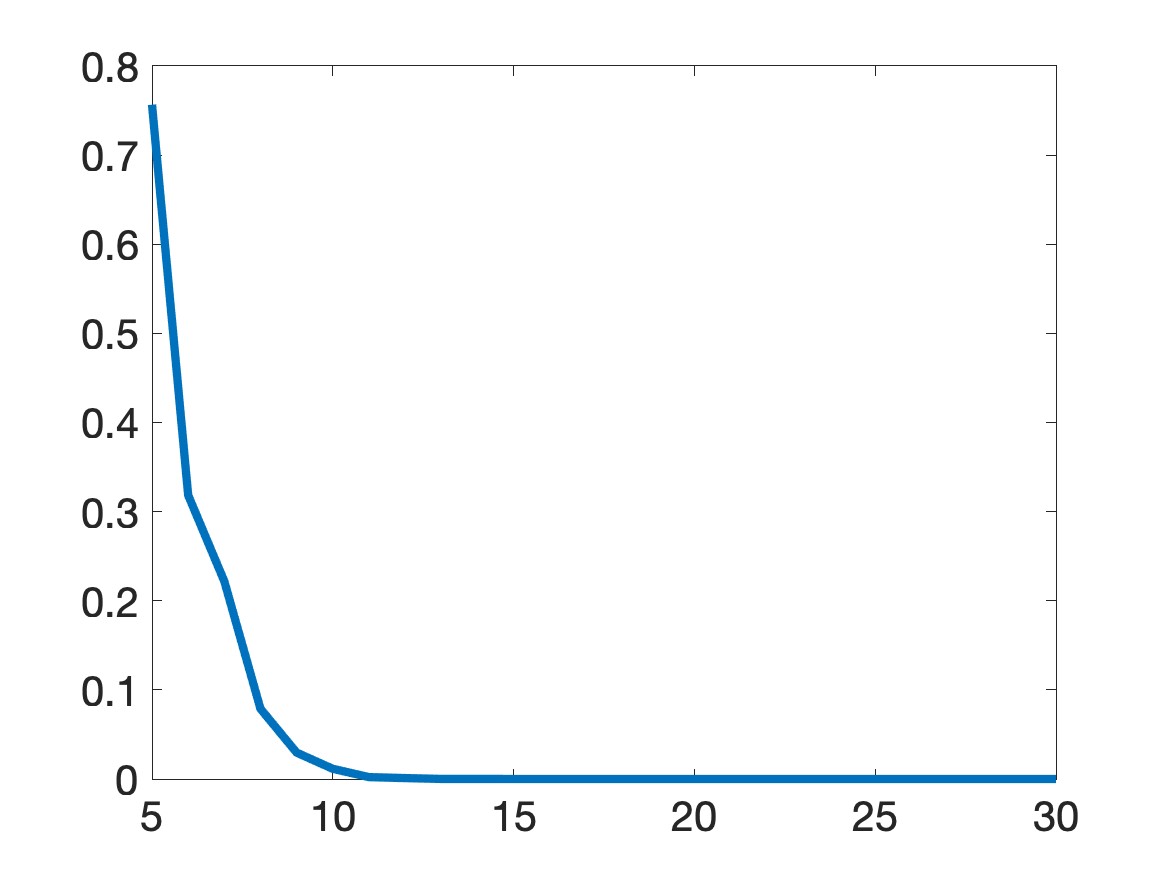}}
	\quad
	\subfloat[]{\includegraphics[width=.3\textwidth]{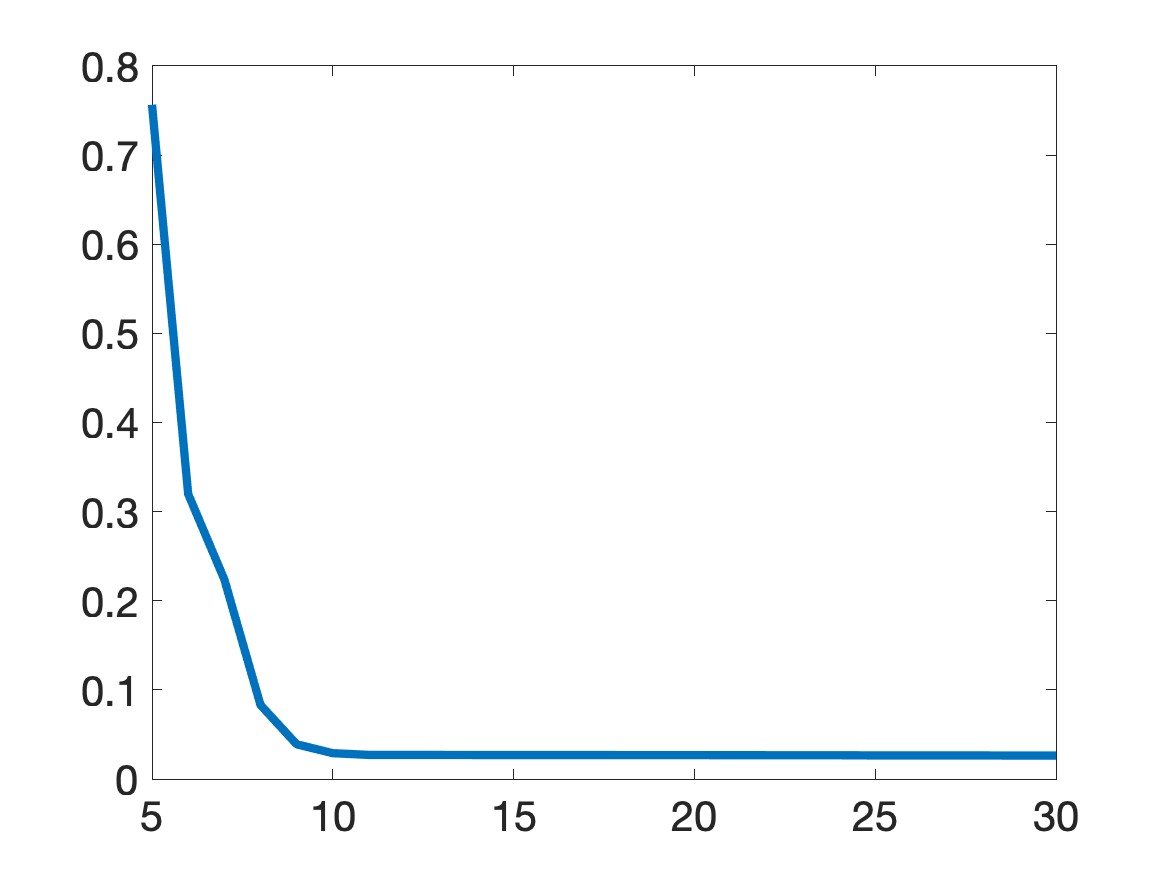}}
	\subfloat[]{\includegraphics[width=.3\textwidth]{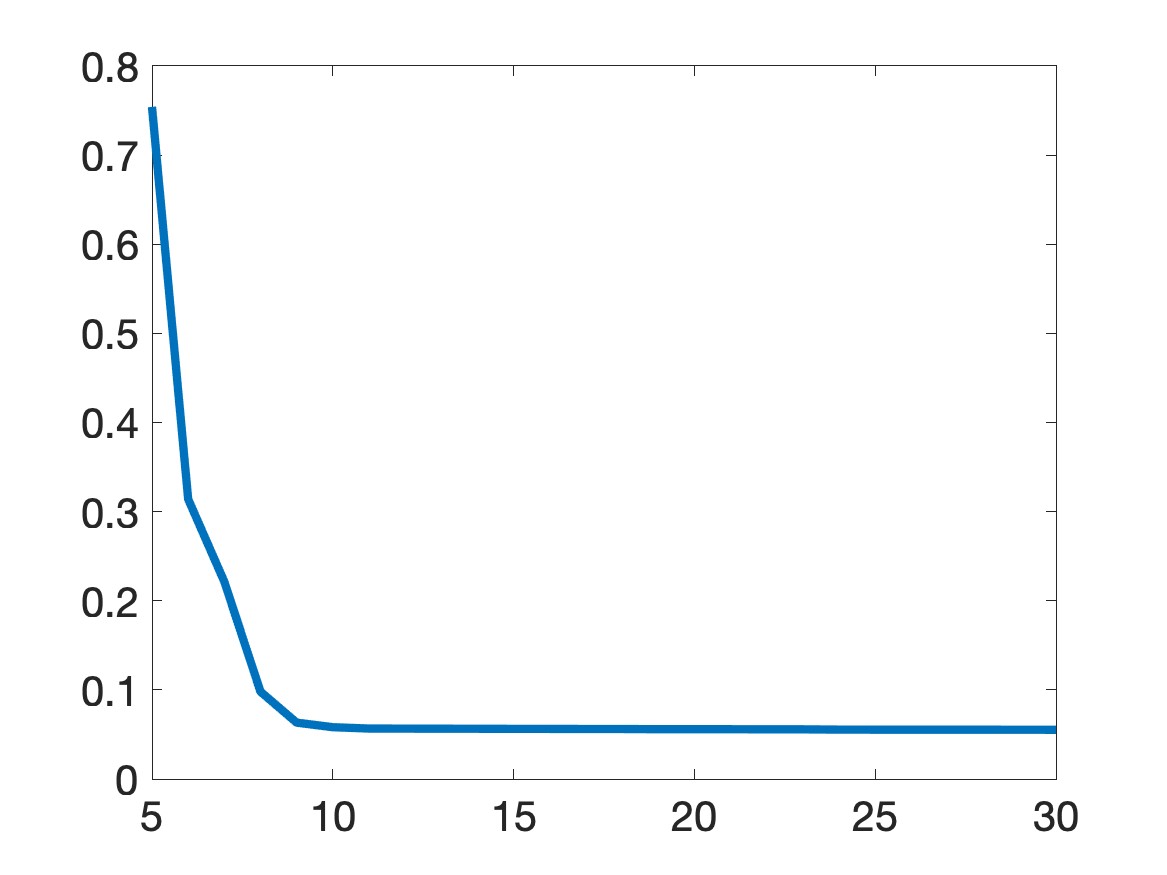}}
	\caption{\label{figN} The graph of the function $\varphi: \{5, \dots, 30\} \to \R$ defined in \eqref{4.2}. The function $\varphi$ is computed from the data $g$ of Test 1 below. (a) $g$ contains $0\%$ of noise; $(b)$ $g$ contains $5\%$ of noise; (c) $g$ contains $10\%$ of noise.}
\end{figure}

We provide here a one-dimensional example of the procedure above to choose the cut-off number for Test 1 below. In this test, since the dimension $d = 1$, $\bN$ becomes the set $\{n_t \geq 1\}$. We only need to determine the cut-off number $N_t$. The graph of $\varphi$ is displayed in Figure \ref{figN} when the data $g$ of Test 1 below is interrupted by $0\%$, $5\%$, and $10\%$ of noise.

It is interesting seeing that the graphs of $\varphi$ looks like an $L-$ curve. 
This observation somewhat confirms that $N_t$ serves as a regularization factor.
One can choose $N_t = 10$ for Test 1 which corresponds to the corner of the $L$.
We experience that choosing higher values of $N_t$, say $10 \leq N_t \leq 15$ delivers similar numerical results. 
However, if we continue to increase $N_t$, say $N_t > 30$, our method breaks out. 
We choose $N_t = 15$ for Test 1.
The procedure to choose the cut-off numbers for other tests is similar.

\begin{Remark}
	It is worth mentioning that it follows from Figure \ref{figN} that the choice of $N_t$ does not depend on the noise levels. For all three noise levels $0\%,$ $5\%$, and $10\%$, the corner of the $L-$curve occurs at   $N_t = 10$. 
\end{Remark}

\subsection{Step \ref{s3} of Algorithm \ref{alg}}

As previously stated, we use the Runge-Kutta method to numerically solve the initial value problem \eqref{2.14}. We omit a detailed explanation of this method in this paper as it can be found in any textbook on ordinary differential equations (ODEs). An alternative approach for computing solutions to \eqref{2.14} is the classical Euler method. However, in order to achieve satisfactory numerical results with the Euler method, it is necessary to choose $N_x$ and $N_t$ to be extremely large, such as $N_x > 2000$ and $N_t > 2000$. Consequently, employing the Runge-Kutta method considerably reduces computational efforts.

In Step \ref{s3}, the Runge-Kutta method is not explicitly implemented. Instead, we utilize the built-in Matlab command "ode45", which is designed for solving a system of ODEs of the form $y'(s) = G(s, y(s))$ subject to an initial condition $y(s_0) = y_0$, where the nonlinear function $G$ and the initial value $y_0$ are known, and $y$ is a column vector-valued function. To adapt our system \eqref{2.14} for use with "ode45", we must first restructure $w_{\bN}$ as a column vector $\mathfrak{w}_N = (\mathfrak{w})_{n = 1}^{N} \in \R^N$, where $N = N_T\prod_{i=2}^d N_i$. The single index $n$ that corresponds to the multi-index $\bn = (n_2, \dots, n_d, n_t) \in \bN$ is given by
\[
	n = (n_t-1)\prod_{i = 2}^d N_i + (n_d-1) \prod_{i = 2}^{d - 1} N_i + \dots + n_2.
\]
The vector $\mathfrak{w}_N$ is referred to as a "line up" version of $w_\bn$.
One can derive a system of ODEs for $\mathfrak{w}_N$ corresponding to \eqref{2.14} easily. This brings us into the context of "ode45". The output of "ode45" for this system of ODEs is denoted by $\mathfrak{w}_N^{\rm comp}$, and the corresponding solution $w_\bN^{\rm comp}$ to \eqref{2.14} can be obtained directly.

\section{Numerical examples} \label{sec 5}

We show numerical results obtained by Algorithm \ref{alg} when $d = 1$ and $d = 2$.

\subsection{The one-dimensional case}

We show two numerical results for the case $d =1$ obtained by Algorithm \ref{alg}.
Recall that in our implementation for the 1D case, we set $\Omega = (0, 1)$ and $T = 1.5$.
The measurements in this section are taken place at $\Gamma_T = \{0\} \times [0, T]$.

\subsubsection{Test 1}\label{sectest1}
 
We test our numerical method for the linear case; i.e., $F\left(x,t,u(x,t),\nabla u(x,t)\right) = 0$ for all $(x, t) \in \Omega_T$.
That means the governing equation is of the form
\begin{equation}
u_t(x, t) = u_{xx}(x, t) \quad \mbox{ for } (x,t) \in \Omega \times (0, T)
\label{eqn test 1}
\end{equation}
The given Cauchy data for this test is
\begin{equation}
	g(t) =  2 \sin(8 t)
	\quad
	\mbox{and }
	\quad
	q(t) = 0
\label{data_test_1}
\end{equation}
for all $t \in (0, T).$
The true solution to Problem \ref{IHCP} in this test is
\begin{equation}
u^{\rm true}(x,t) = e^{2x}\sin\left(8t+2x\right)+e^{-2x}\sin\left(8 t-2x\right).
\end{equation}

\begin{figure}[h!]
	\subfloat[The function $u_{\rm true}$]{\includegraphics[width = .3\textwidth]{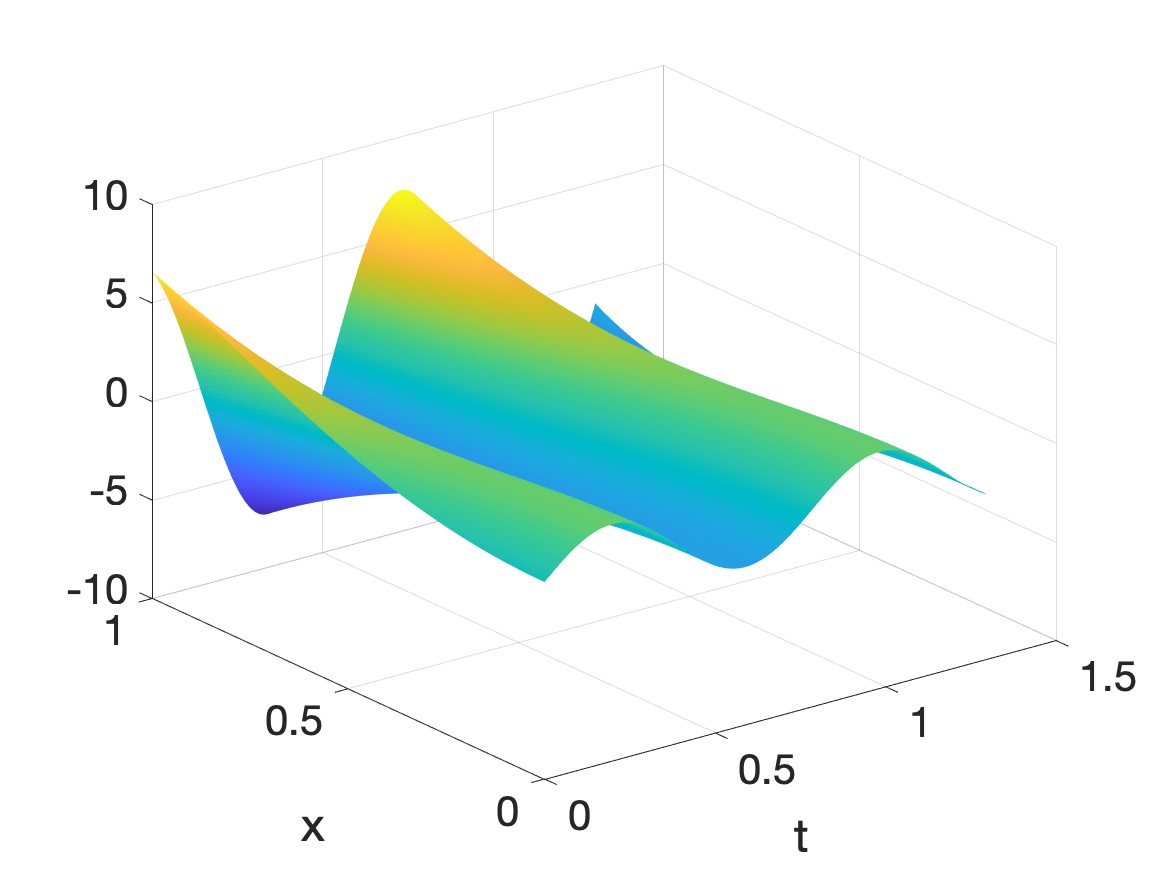}}
	
	\subfloat[The function $u_{\rm comp}$, $\delta = 0\%$]{\includegraphics[width = .3\textwidth]{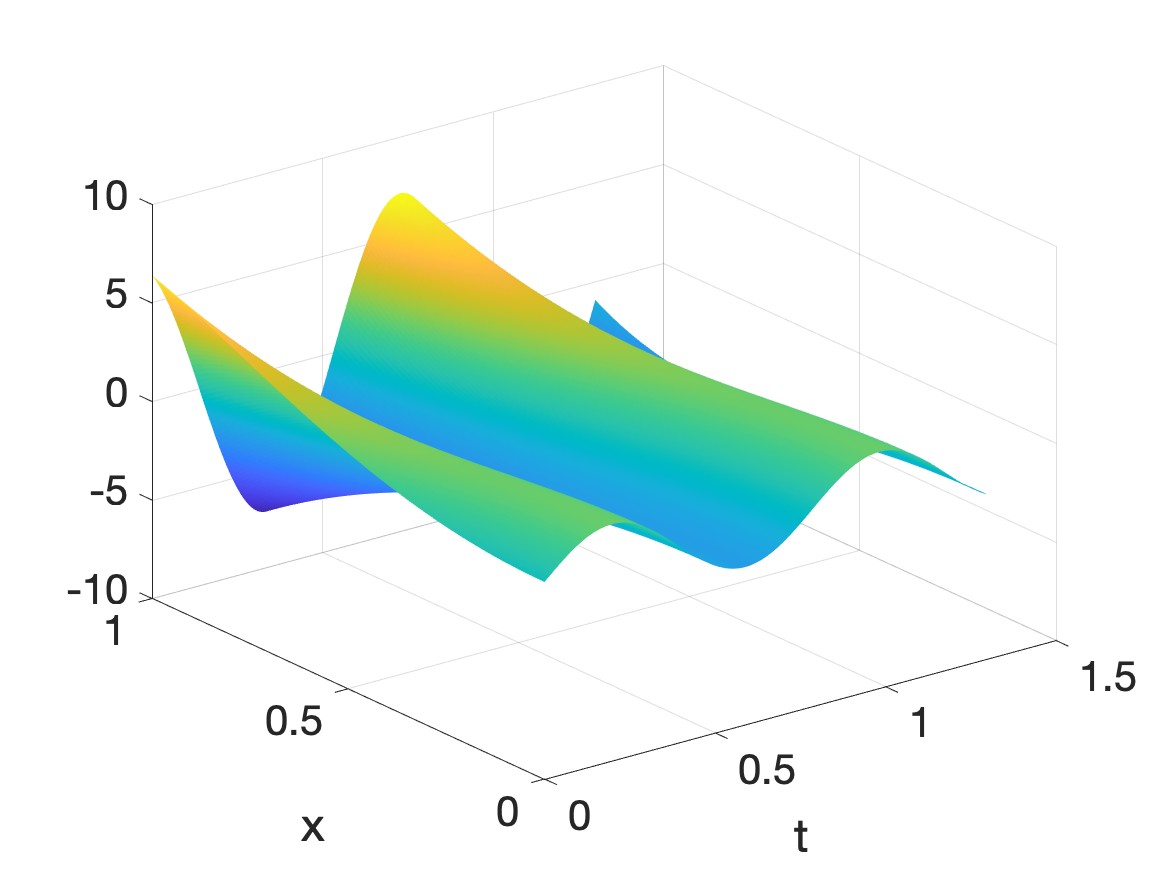}}
	\quad
	\subfloat[The function $u_{\rm comp}$, $\delta = 5\%$]{\includegraphics[width = .3\textwidth]{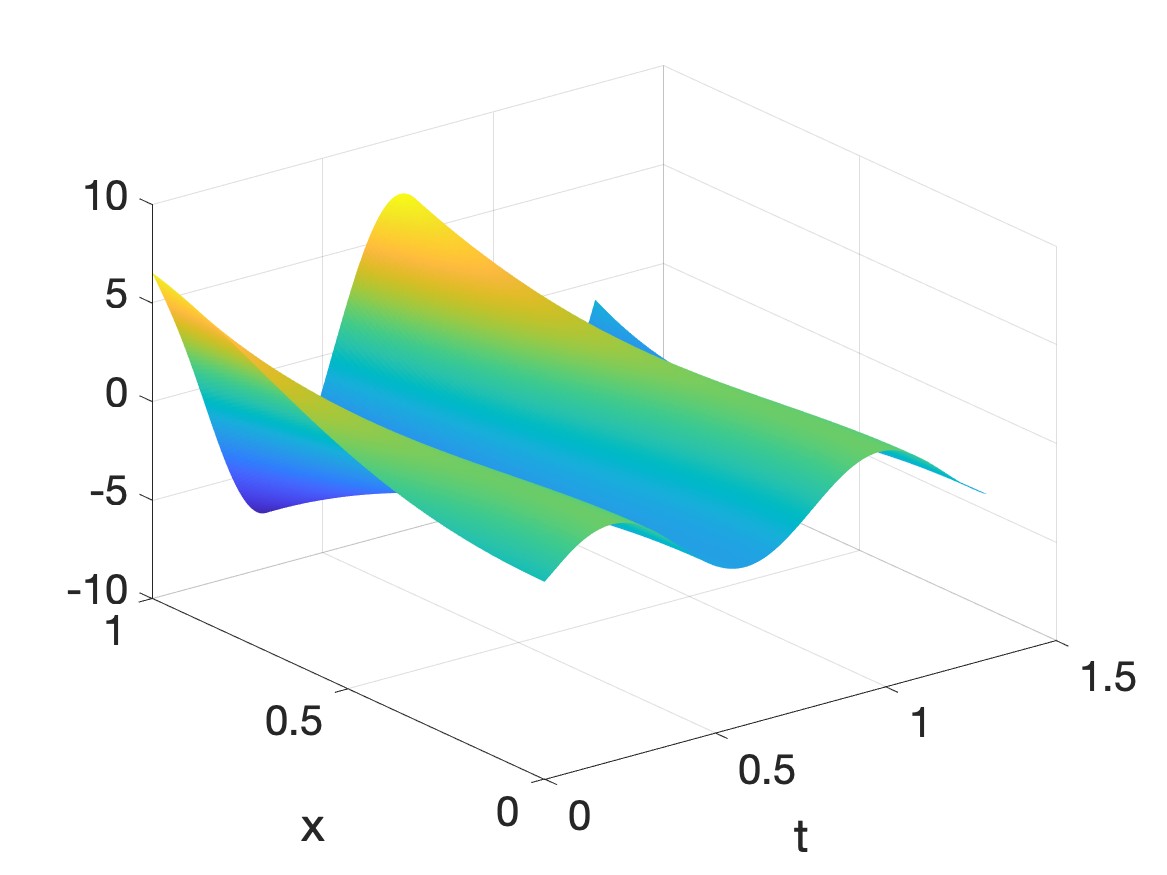}}
	\quad
	\subfloat[The function $u_{\rm comp}$, $\delta = 10\%$]{\includegraphics[width = .3\textwidth]{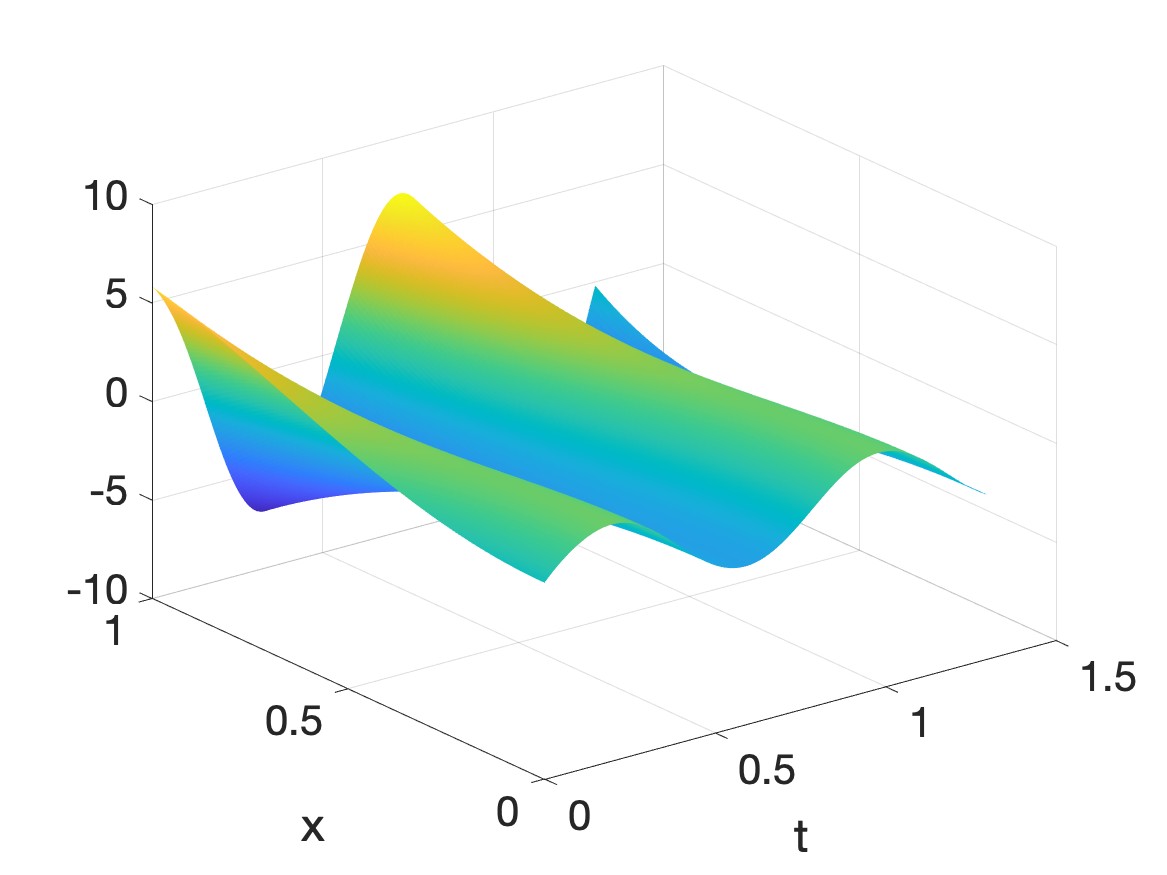}}
	
	\subfloat[The function $\frac{|u_{\rm true} - u_{\rm comp}|}{\|u_{\rm true}\|_{L^\infty}}$, $\delta = 0\%$]{\includegraphics[width = .3\textwidth]{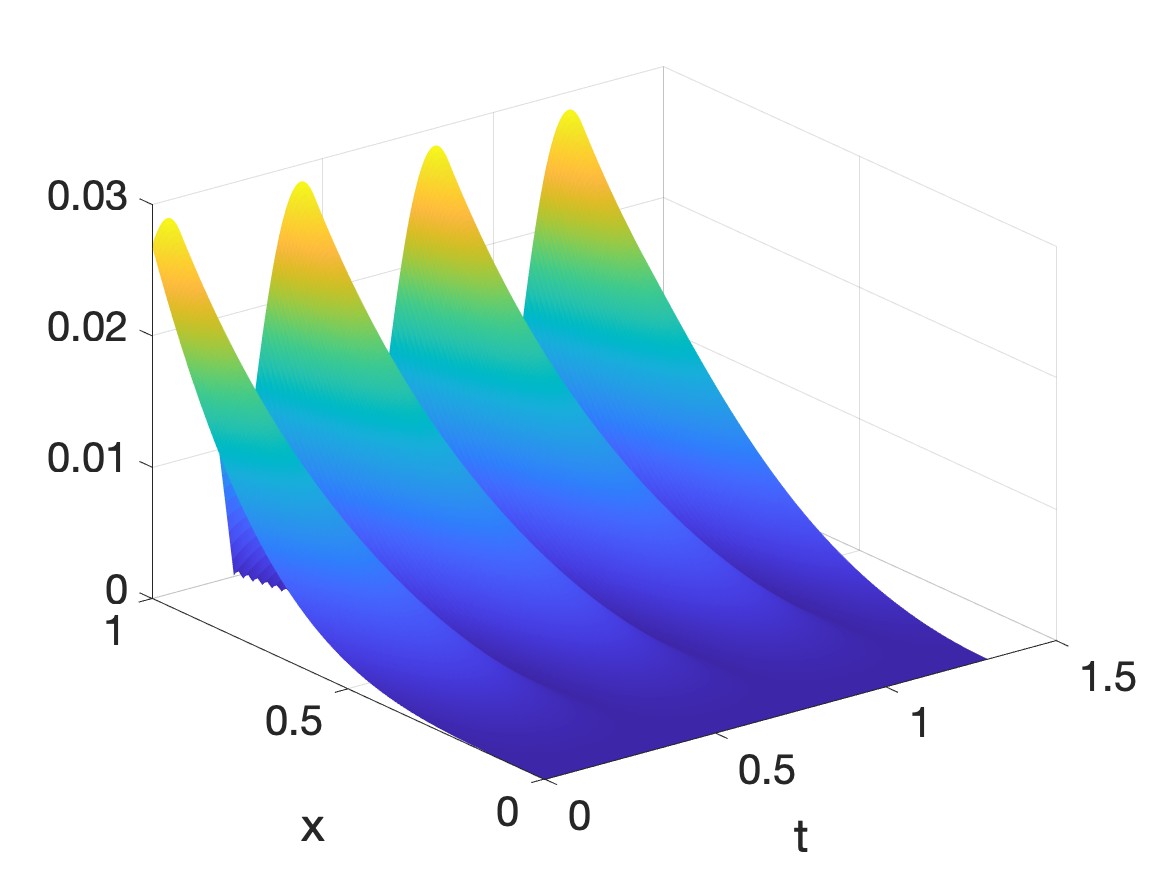}}
	\quad
	\subfloat[The function $\frac{|u_{\rm true} - u_{\rm comp}|}{\|u_{\rm true}\|_{L^\infty}}$, $\delta = 5\%$]{\includegraphics[width = .3\textwidth]{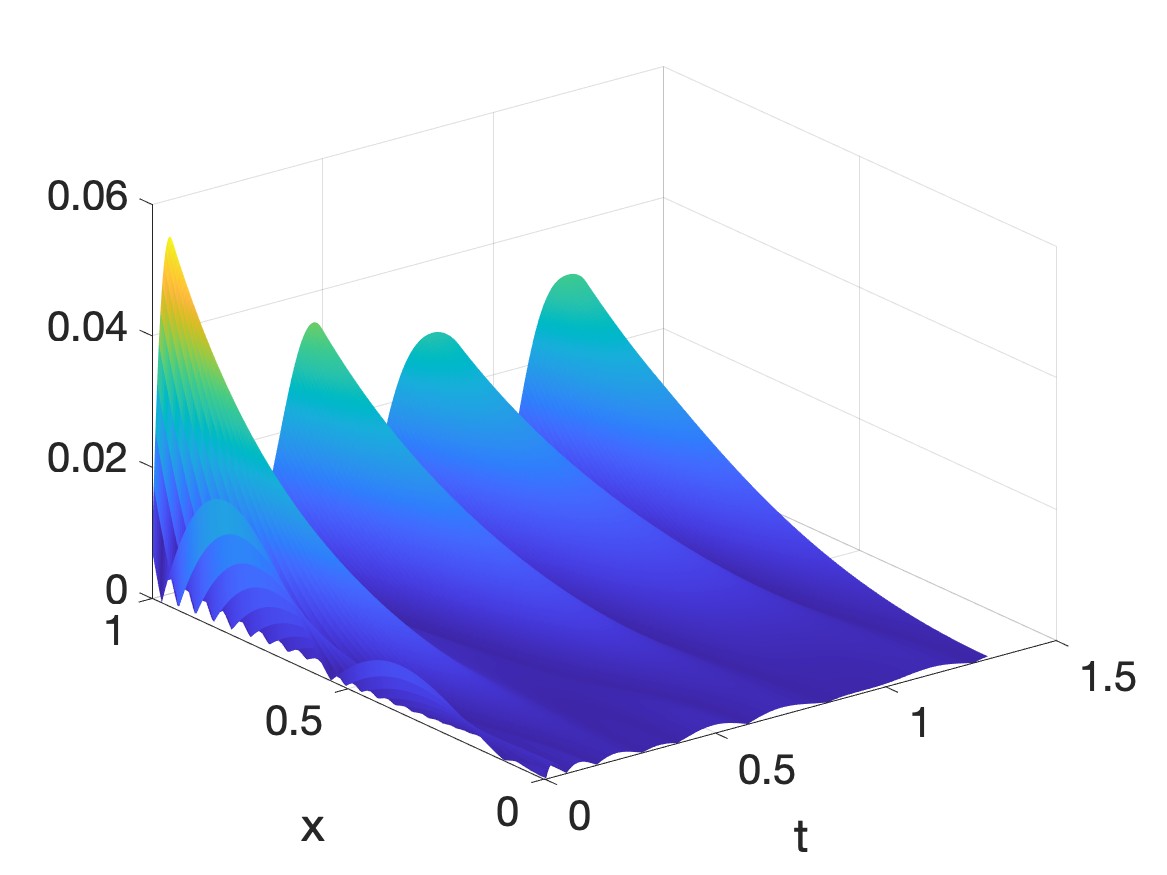}}
	\quad
	\subfloat[The function $\frac{|u_{\rm true} - u_{\rm comp}|}{\|u_{\rm true}\|_{L^\infty}}$, $\delta = 10\%$]{\includegraphics[width = .3\textwidth]{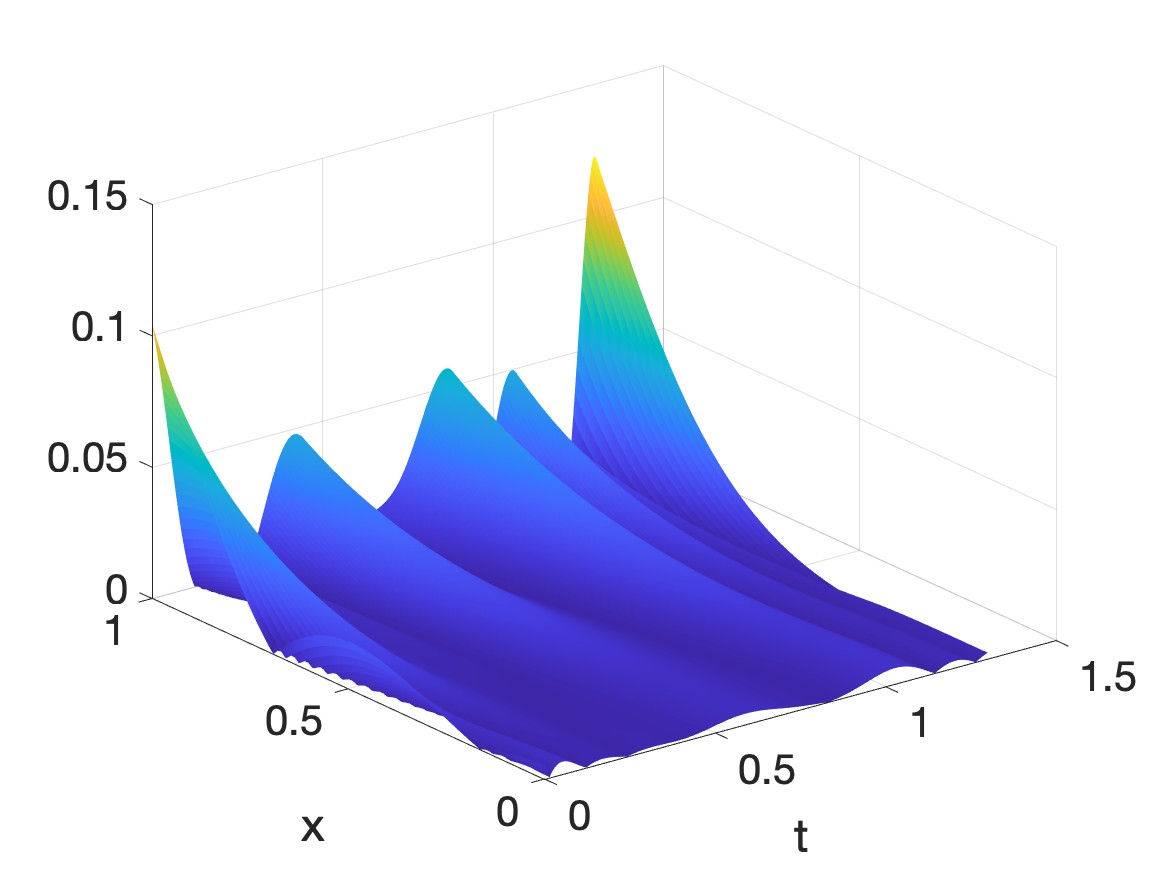}}
	\caption{\label{figtest1}Numerical results of Test 1 when the noise levels are $0,$ $5\%$, and $10\%$. The computation error primarily arises at $x=1$, which is away from the location of the measurement point, $x=0$.
	}
\end{figure}

Figure \ref{figtest1} displays the true and computed solutions, along with their differences, for Test 1. 

Algorithm \ref{alg} successfully provides out-of-expectation computed solutions to Problem \ref{IHCP} for the linear parabolic equation \eqref{eqn test 1} with the Cauchy data given in \eqref{data_test_1} containing random noise. 
The relative errors $\frac{\|u_{\rm true} - u_{\rm comp}\|_{L^2}}{\|u_{\rm true}\|_{L^2}}$ on  $(0, 1) \times (0, 1.3)$ are $0.32\%$, $0.50\%$, and $1.59\%$ when the noise level $\delta$  is $0\%$, $5\%$, and $10\%$ respectively.
It is interesting to note that the computation errors are compatible with the noise levels.

\subsubsection{Test 2}
\begin{figure}[h!]
	\subfloat[The function $u_{\rm true}$]{\includegraphics[width = .3\textwidth]{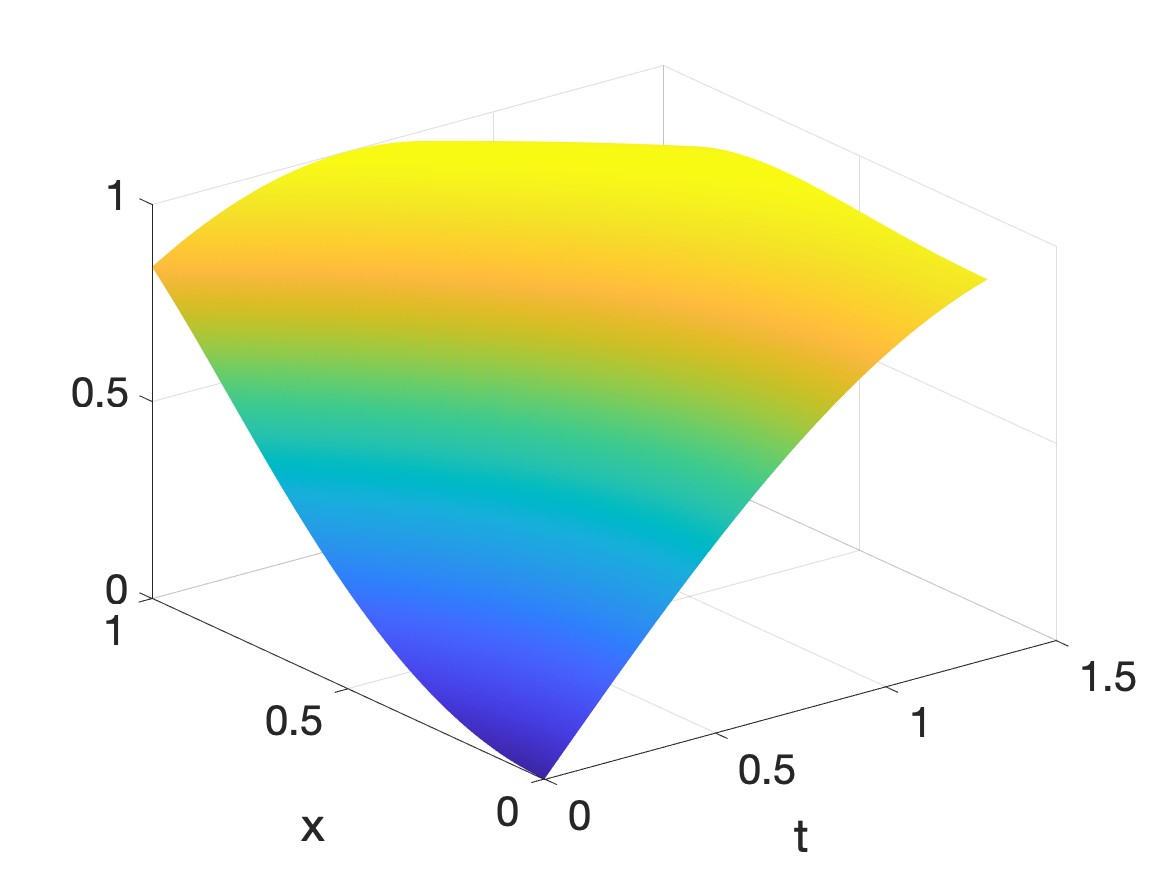}}
	
	\subfloat[The function $u_{\rm comp}$, $\delta = 0\%$]{\includegraphics[width = .3\textwidth]{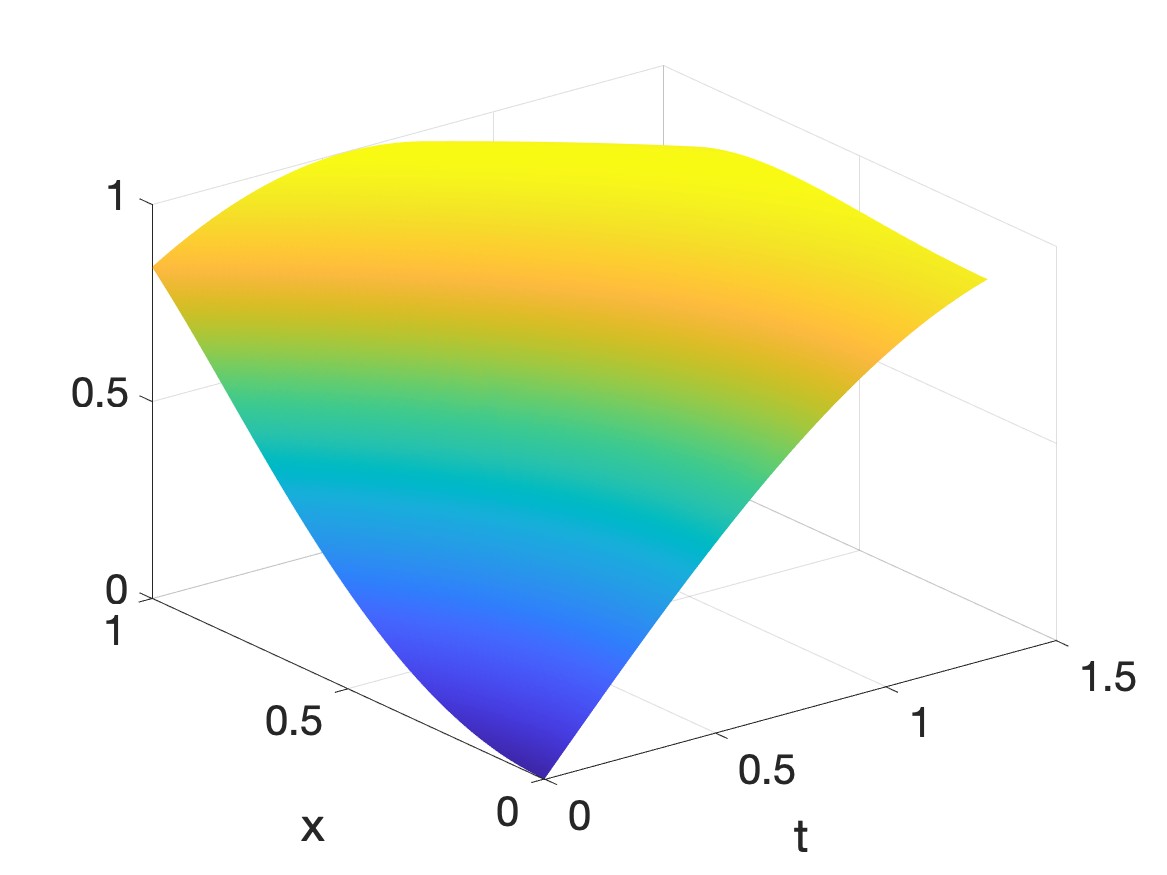}}
	\quad
	\subfloat[The function $u_{\rm comp}$, $\delta = 5\%$]{\includegraphics[width = .3\textwidth]{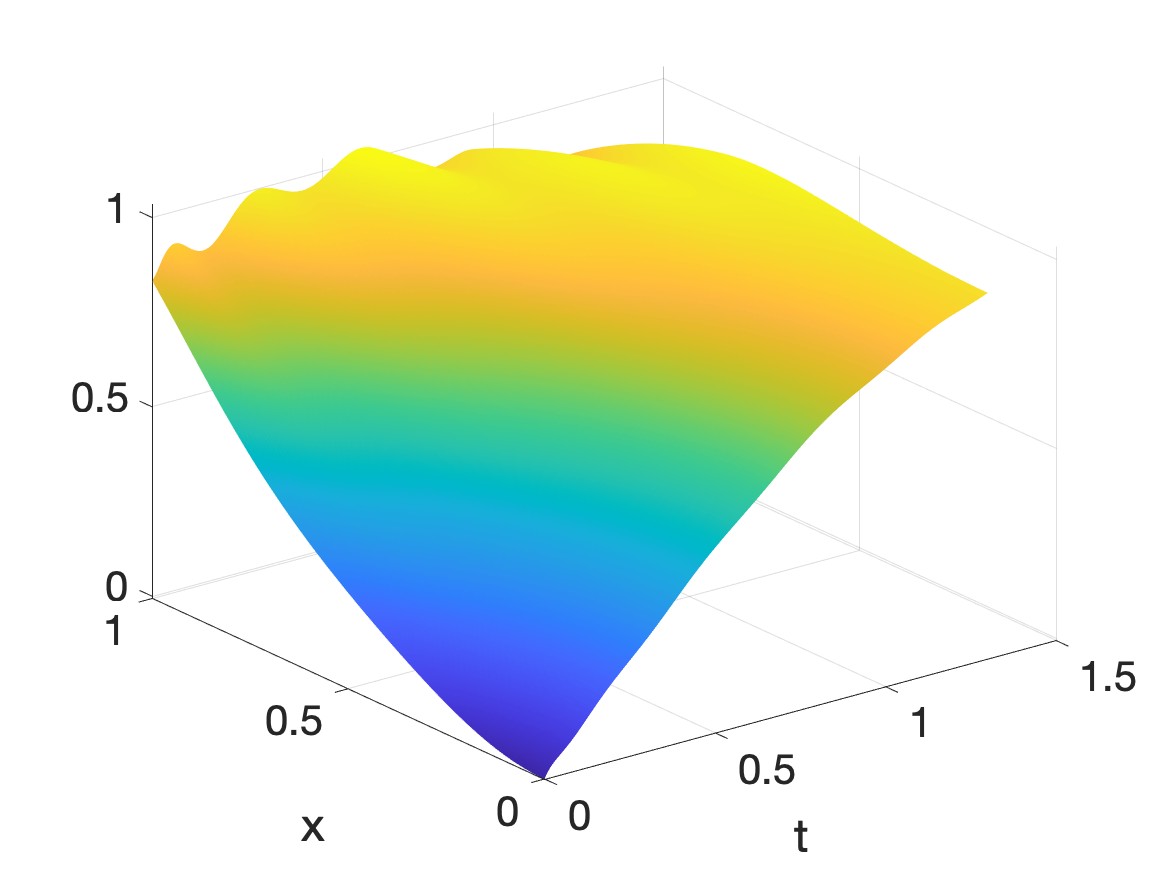}}
	\quad
	\subfloat[The function $u_{\rm comp}$, $\delta = 10\%$]{\includegraphics[width = .3\textwidth]{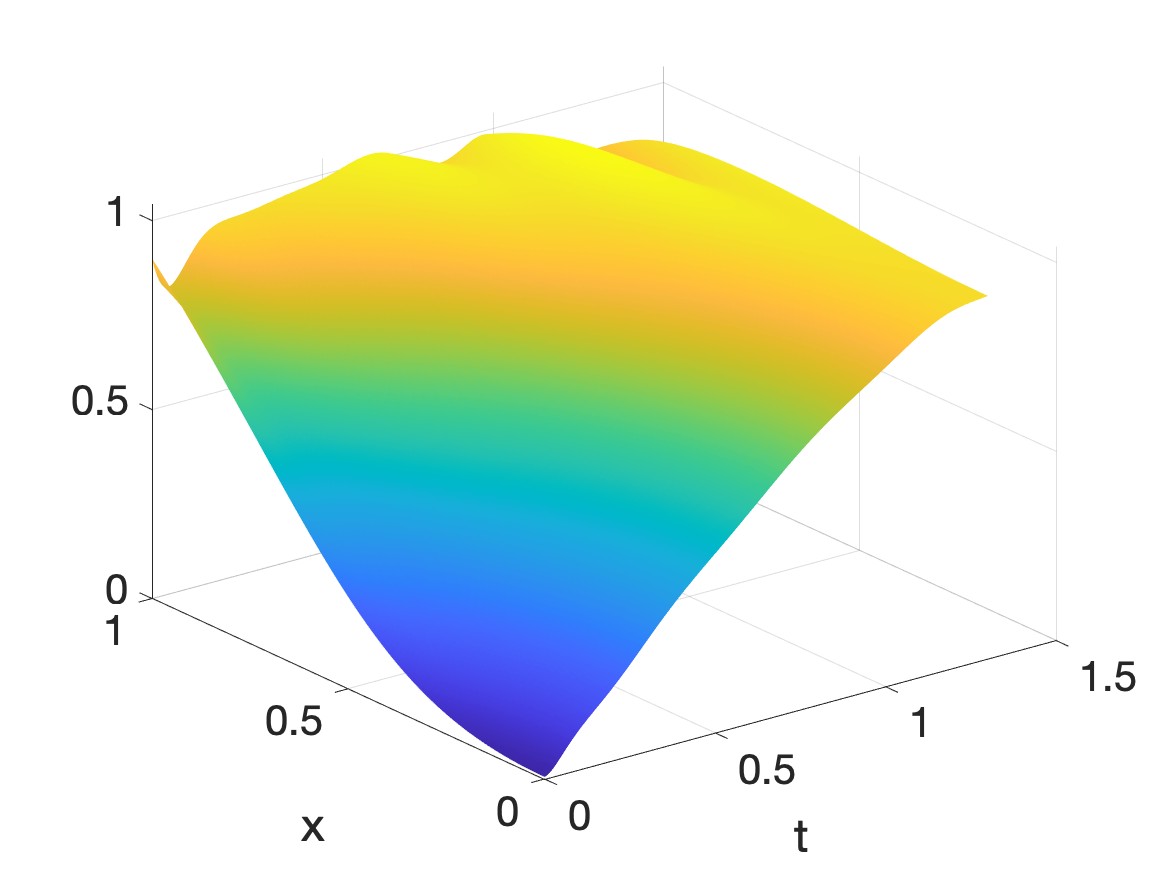}}
	
	\subfloat[The function $\frac{|u_{\rm true} - u_{\rm comp}|}{\|u_{\rm true}\|_{L^\infty}}$, $\delta = 0\%$]{\includegraphics[width = .3\textwidth]{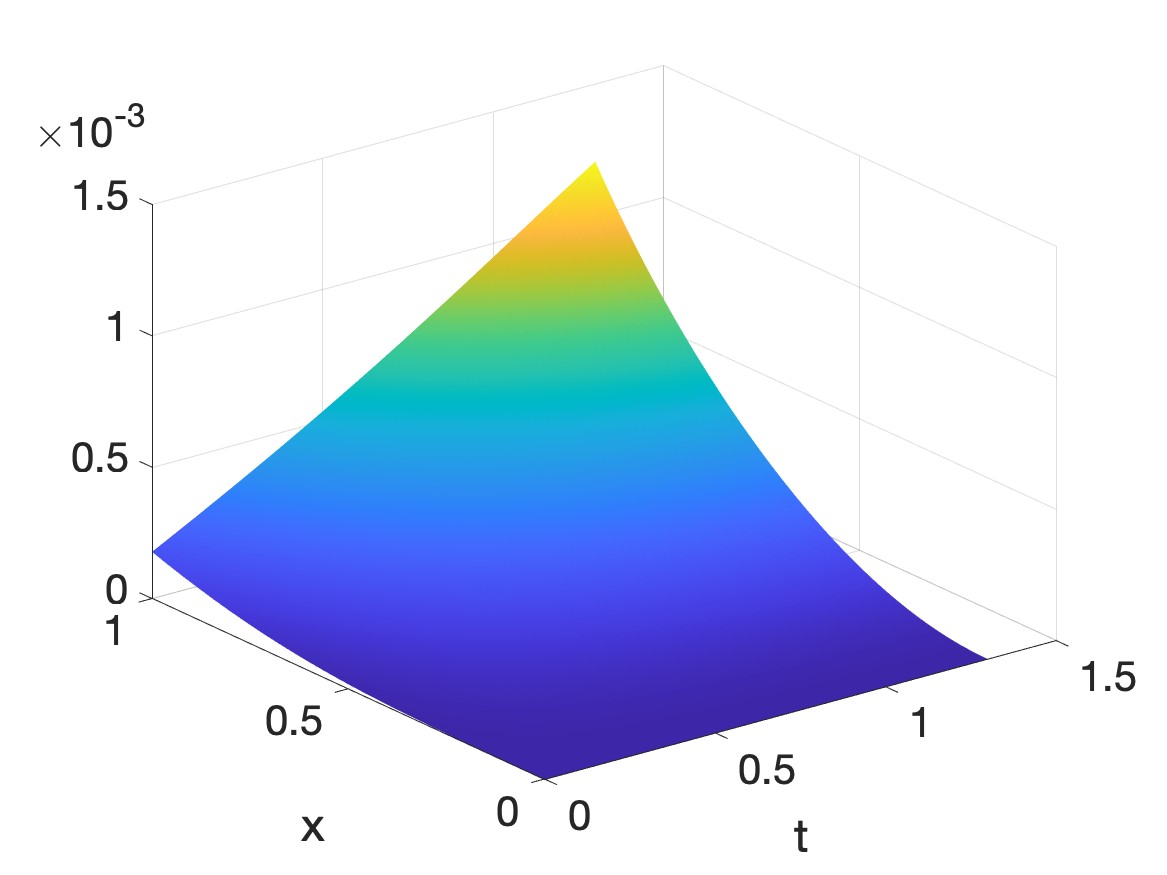}}
	\quad
	\subfloat[The function $\frac{|u_{\rm true} - u_{\rm comp}|}{\|u_{\rm true}\|_{L^\infty}}$, $\delta = 5\%$]{\includegraphics[width = .3\textwidth]{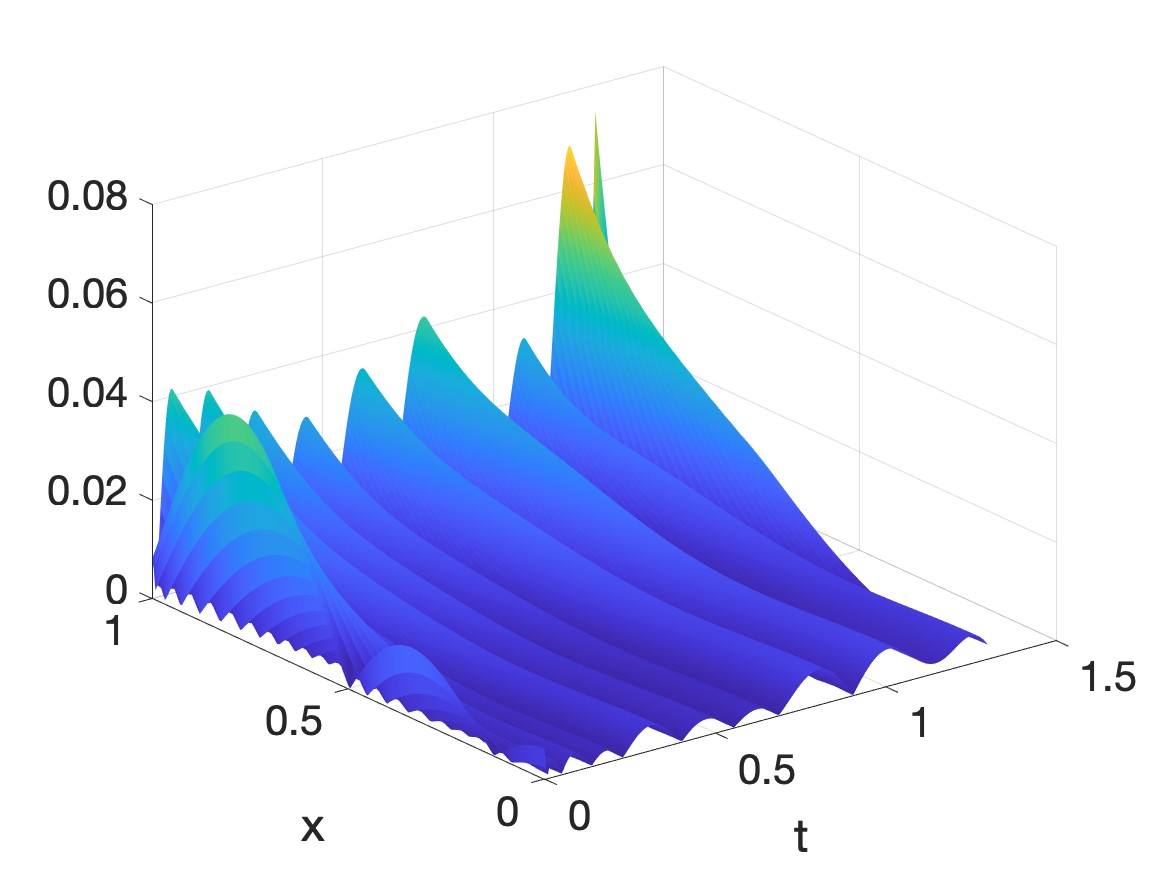}}
	\quad
	\subfloat[The function $\frac{|u_{\rm true} - u_{\rm comp}|}{\|u_{\rm true}\|_{L^\infty}}$, $\delta = 10\%$]{\includegraphics[width = .3\textwidth]{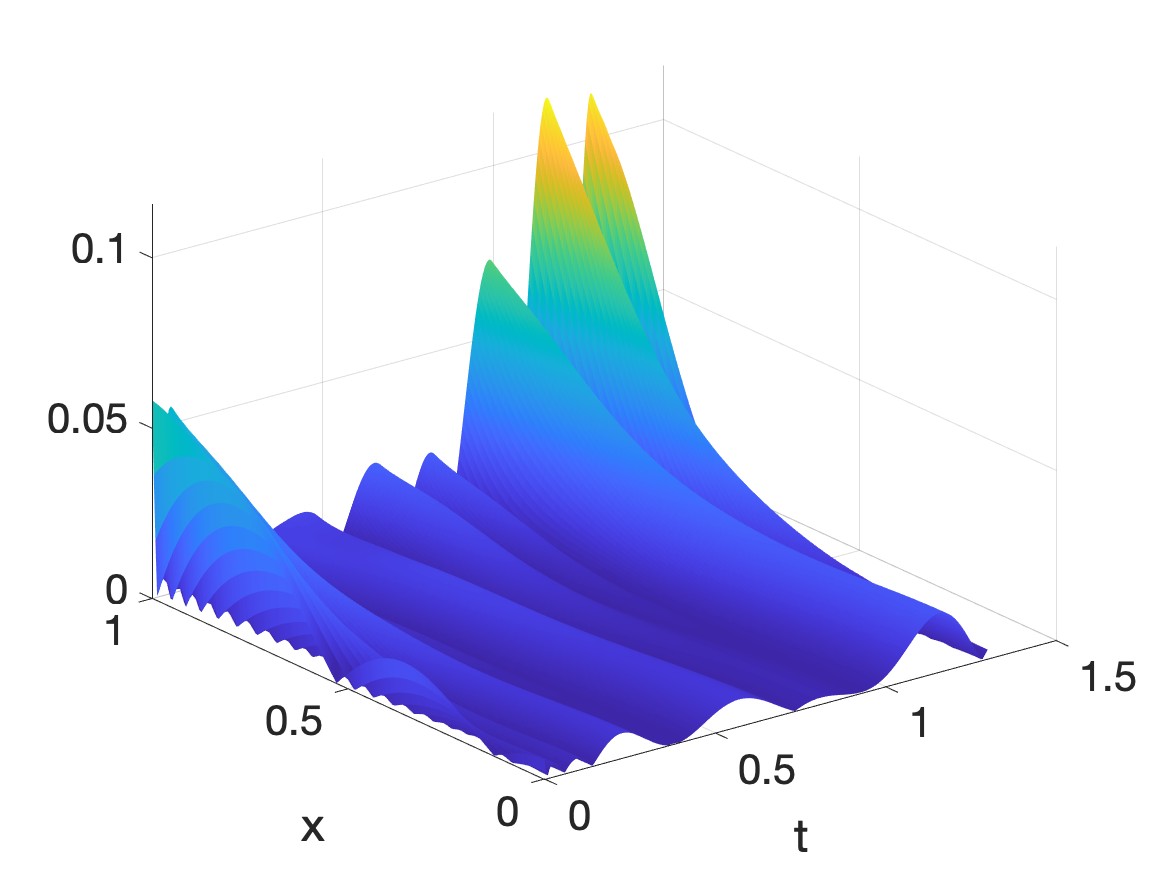}}
	\caption{\label{figtest2}Numerical results of Test 2 when the noise levels are $0,$ $5\%$, and $10\%$.
	The computation error primarily arises at $x=1$, which is away from the location of the measurement point, $x=0$.
	}
\end{figure}

We test the performance of Algorithm \ref{alg} in the nonlinear case.
The nonlinear parabolic equation for this test is
\begin{equation}
	u_t(x, t) = u_{xx}(x, t) + |u_x(x, t)|^2 + 4x^2 u(x, t) - 4x^2 \cos^2(x^2 + t) - \cos(x^2 + t)
	\label{5.4}
\end{equation}
in $\Omega_T$ and the given data are
\begin{equation}
	g(t) = \sin(t)
	\quad
	\mbox{and}
	\quad
	q(t) = 0
	\label{5.5}
\end{equation}
for $t \in (0, T)$.
The true solution in this test is
\begin{equation}
u_{\rm true}(x,t) = \sin(x^2+t) 
\quad
\mbox{ for } (x,t) \in \Omega_T.
\end{equation}
Figure \ref{figtest2} displays the true and computed solutions, along with their differences, for Test 2. 

The nonlinearity $|u_x(x, t)|^2 $ in \eqref{5.4} in this test grows as a quadratic function, which is significantly faster than a linear function. Therefore, numerically solving \eqref{5.4} with given noisy Cauchy data is extremely challenging.
Despite this, our reductional dimension reduction method still delivers satisfactory approximations of the solution.
The relative errors $\frac{\|u_{\rm true} - u_{\rm comp}\|_{L^2}}{\|u_{\rm true}\|_{L^2}}$ on $\Omega \times (0, T') = (0, 1) \times (0, 1.3)$ are $0.00\%$, $0.04\%$, and $0.14\%$ when the noise level $\delta$  is $0\%$, $5\%$, and $10\%$ respectively.
Like in Test 1, computation errors in Test 2 are compatible with the noise levels.

\subsection{The two-dimensional case}
We show two numerical results for the case $d =2$ obtained by Algorithm \ref{alg}.
Recall that in our implementation for the 2D case, we set $\Omega = (0, 0.5) \times (-1, 1)$ and $T = 1.5$.
The measurements in this section are taken place at $\Gamma_T = \{0\} \times (-1, 1) \times (0, T)$.

\subsubsection{Test 3}

We solve the following equation 
\begin{equation}
	u_t(\x, t) = \Delta u(\x, t) + |u(\x, t)|^2 + 4|\x|^2 u(\x, t) - \cos(|\x|^2 + t) + 3 \sin(|\x|^2 + t)
	\label{5.7}
\end{equation}
for $\x = (x, y) \in \Omega$ and $t \in (0, T)$.
The given data on $\Gamma_T$ are
\begin{equation}
	g(y, t) = \cos(y^2 + t)
	\quad
	\mbox{and}
	\quad
	q(y, t) = 0
	\label{5.8}
\end{equation}
for all $(y, t) \in (-1, 1) \times (0, 1.5).$
The true solution is 
\begin{equation}
	u_{\rm true} = \cos(x^2 + y^2 + t).
\end{equation}

\begin{figure}[h!]
	\subfloat[The function $u_{\rm true}$]{\includegraphics[width = .3\textwidth]{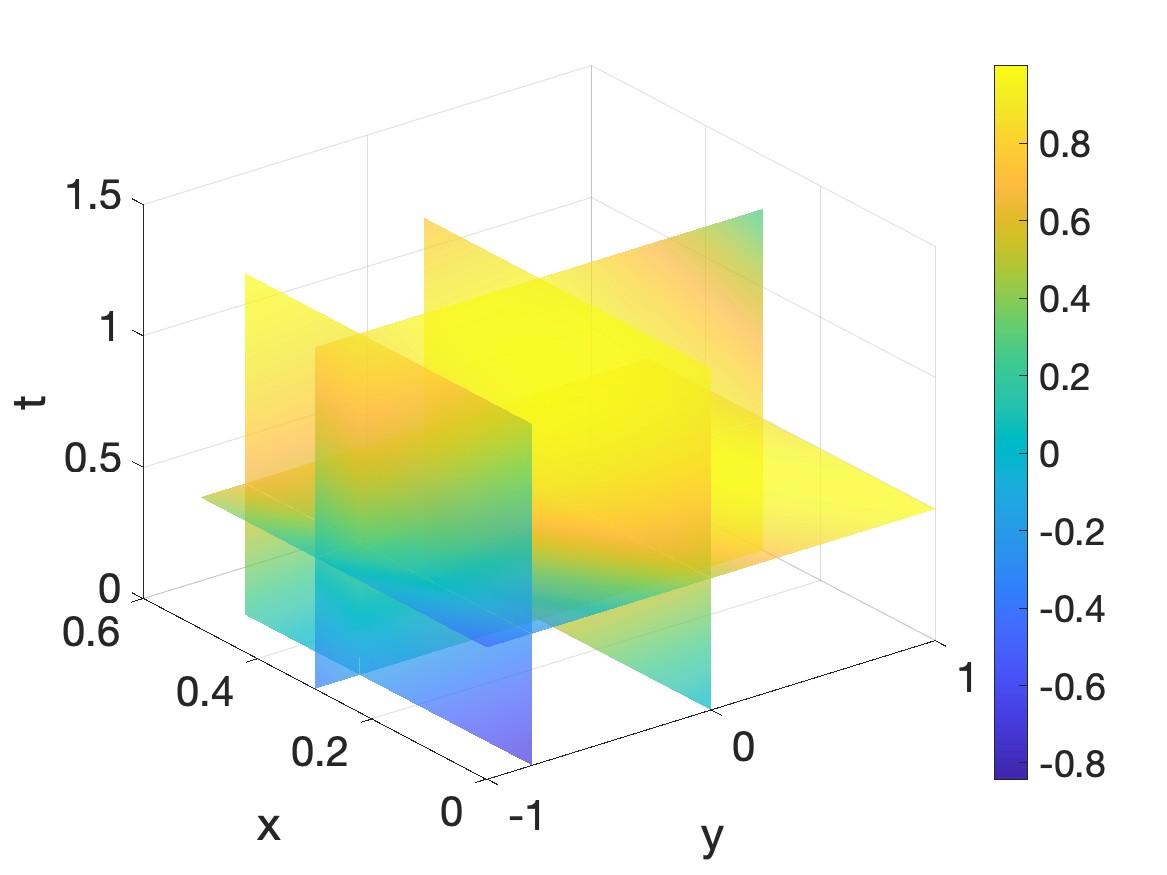}}
	
	\subfloat[The function $u_{\rm comp}$, $\delta = 0\%$]{\includegraphics[width = .3\textwidth]{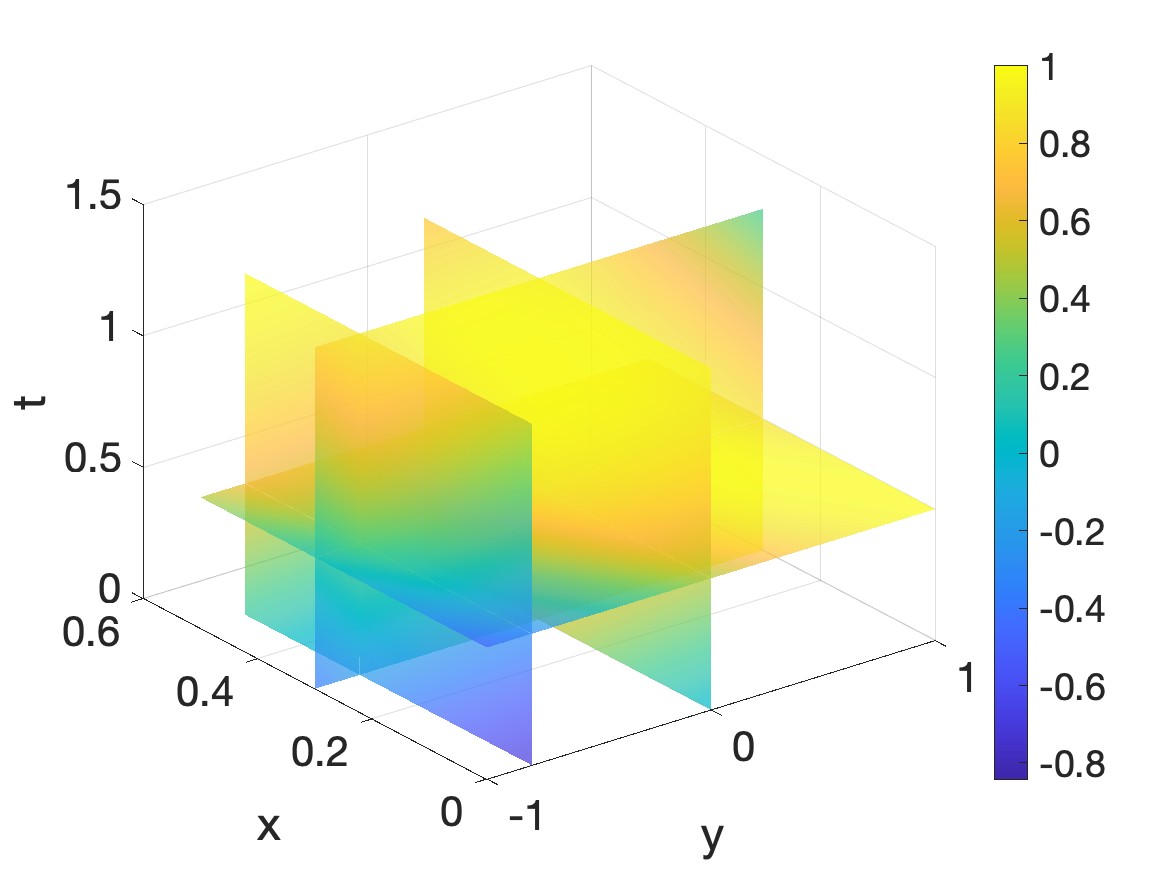}}
	\quad
	\subfloat[The function $u_{\rm comp}$, $\delta = 5\%$]{\includegraphics[width = .3\textwidth]{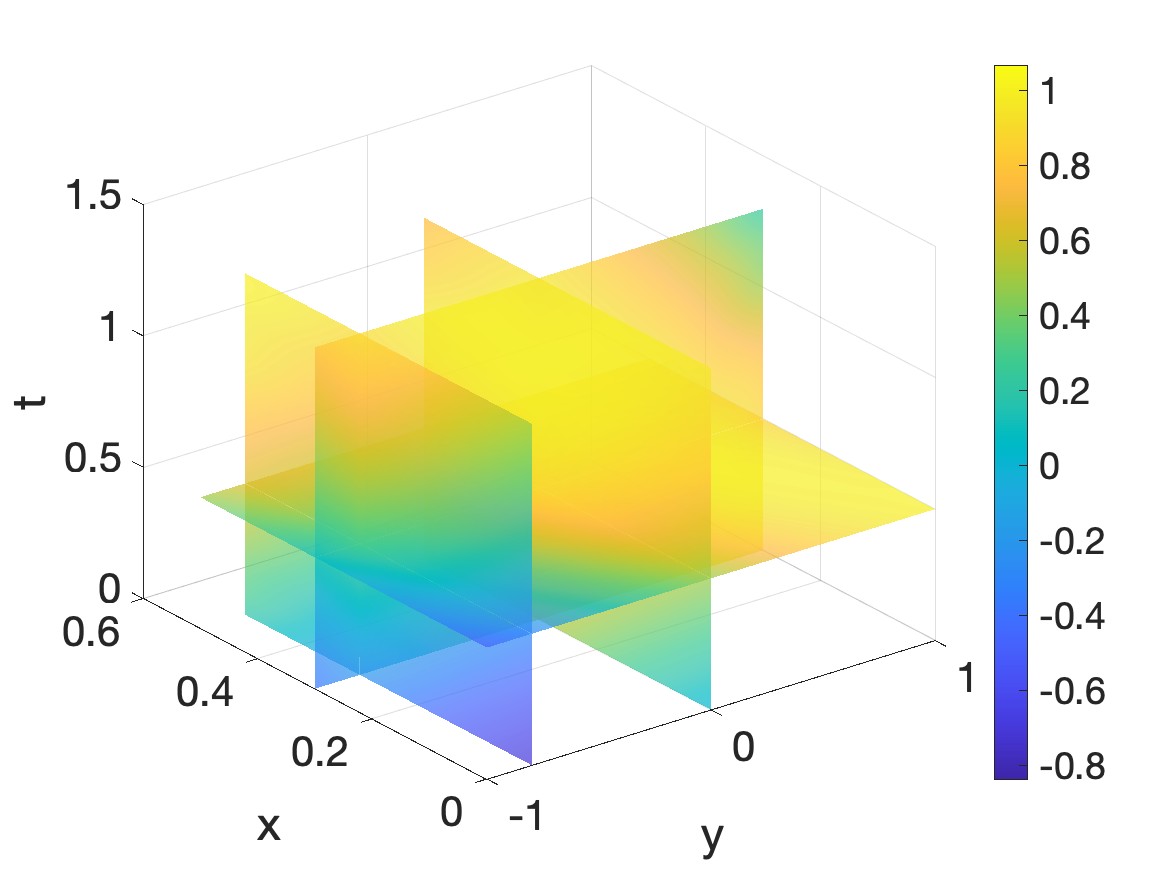}}
	\quad
	\subfloat[The function $u_{\rm comp}$, $\delta = 10\%$]{\includegraphics[width = .3\textwidth]{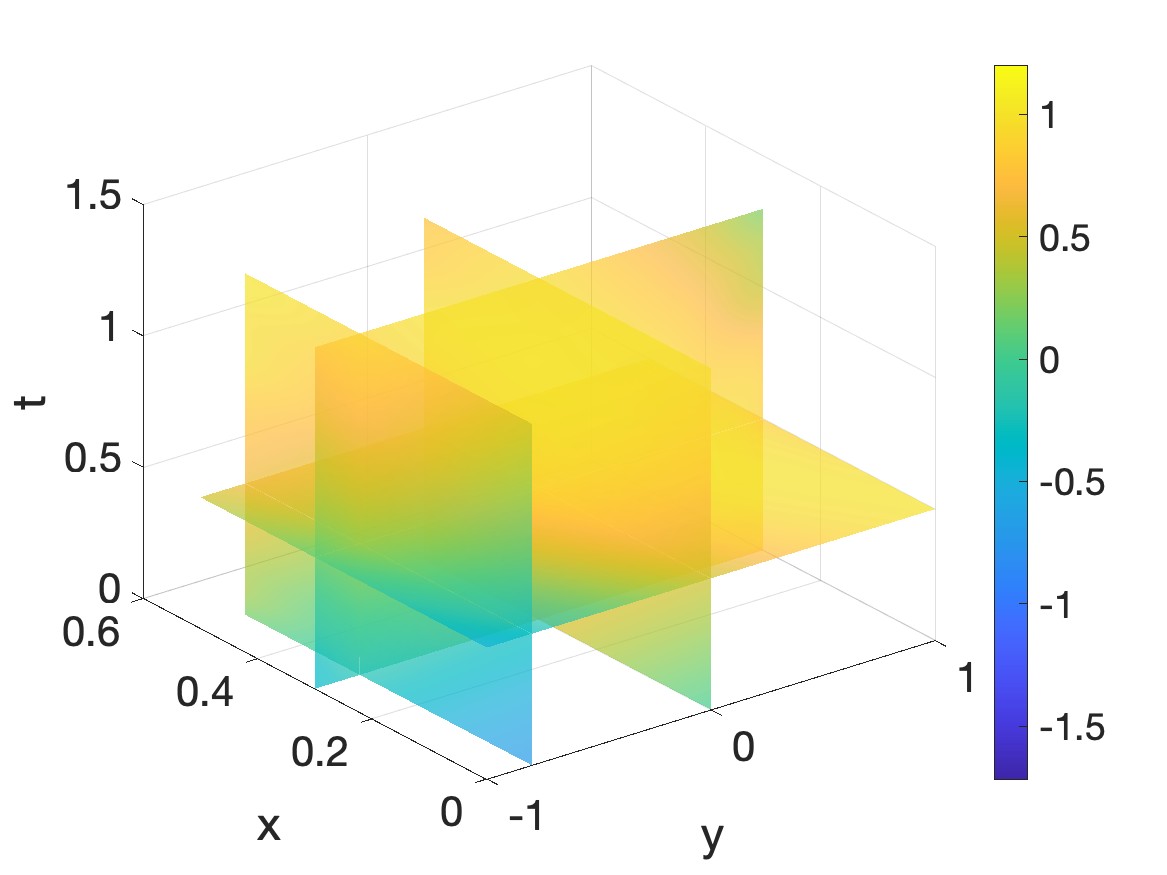}}
	
	\subfloat[The function $\frac{|u_{\rm true} - u_{\rm comp}|}{\|u_{\rm true}\|_{L^\infty}}$, $\delta = 0\%$]{\includegraphics[width = .3\textwidth]{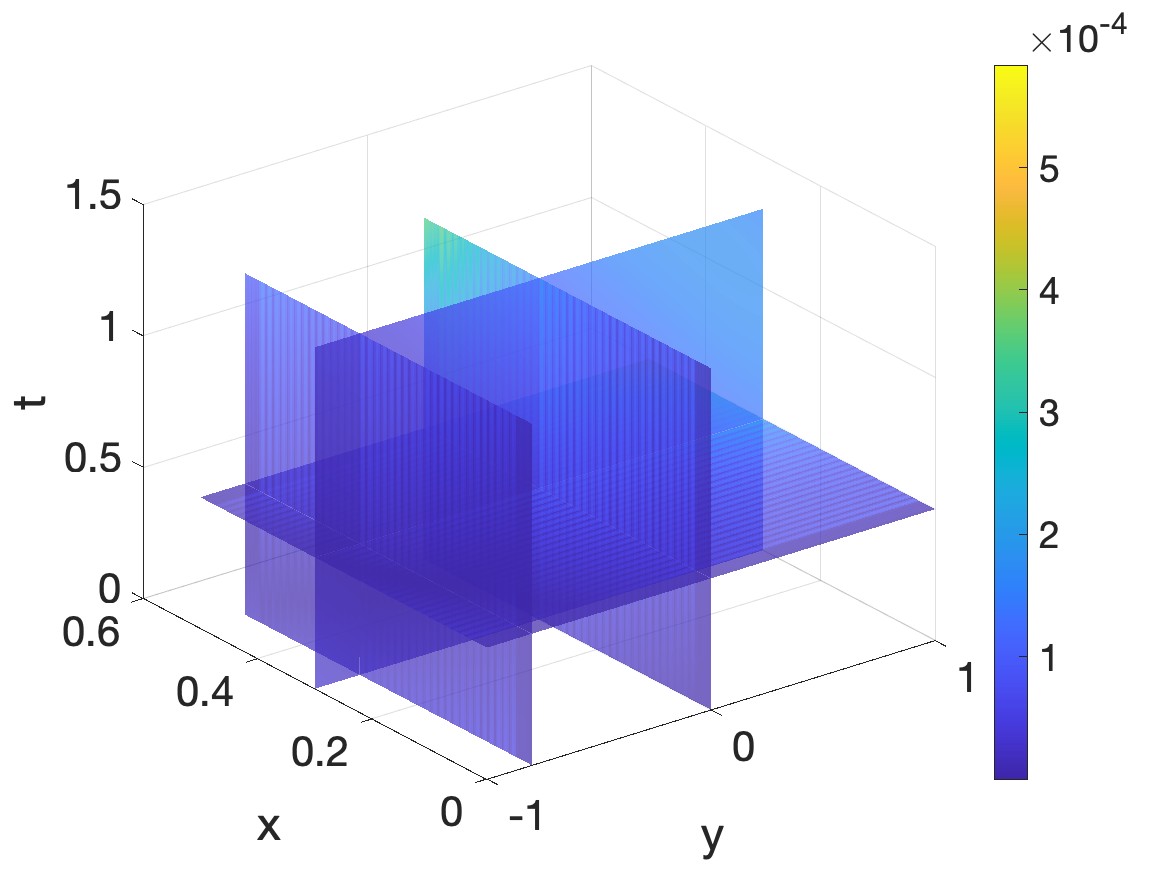}}
	\quad
	\subfloat[The function $\frac{|u_{\rm true} - u_{\rm comp}|}{\|u_{\rm true}\|_{L^\infty}}$, $\delta = 5\%$]{\includegraphics[width = .3\textwidth]{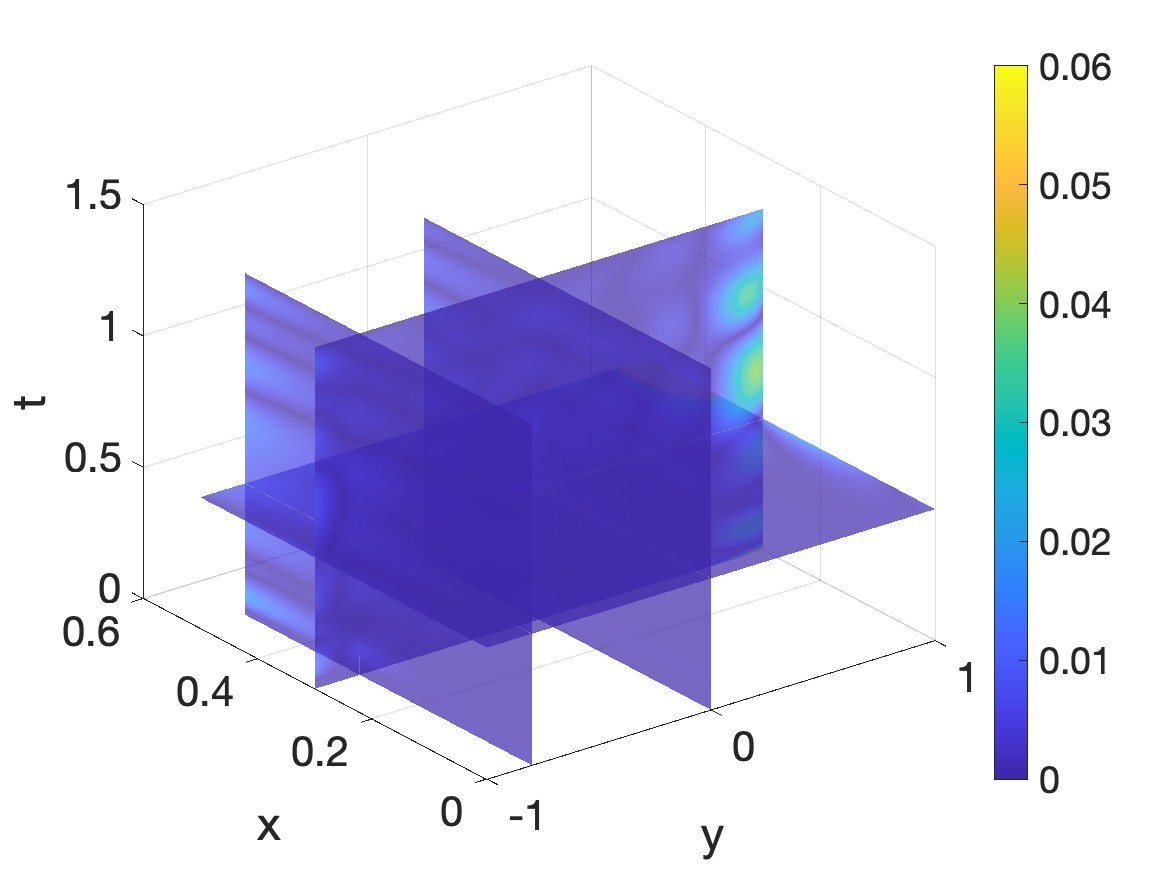}}
	\quad
	\subfloat[The function $\frac{|u_{\rm true} - u_{\rm comp}|}{\|u_{\rm true}\|_{L^\infty}}$, $\delta = 10\%$]{\includegraphics[width = .3\textwidth]{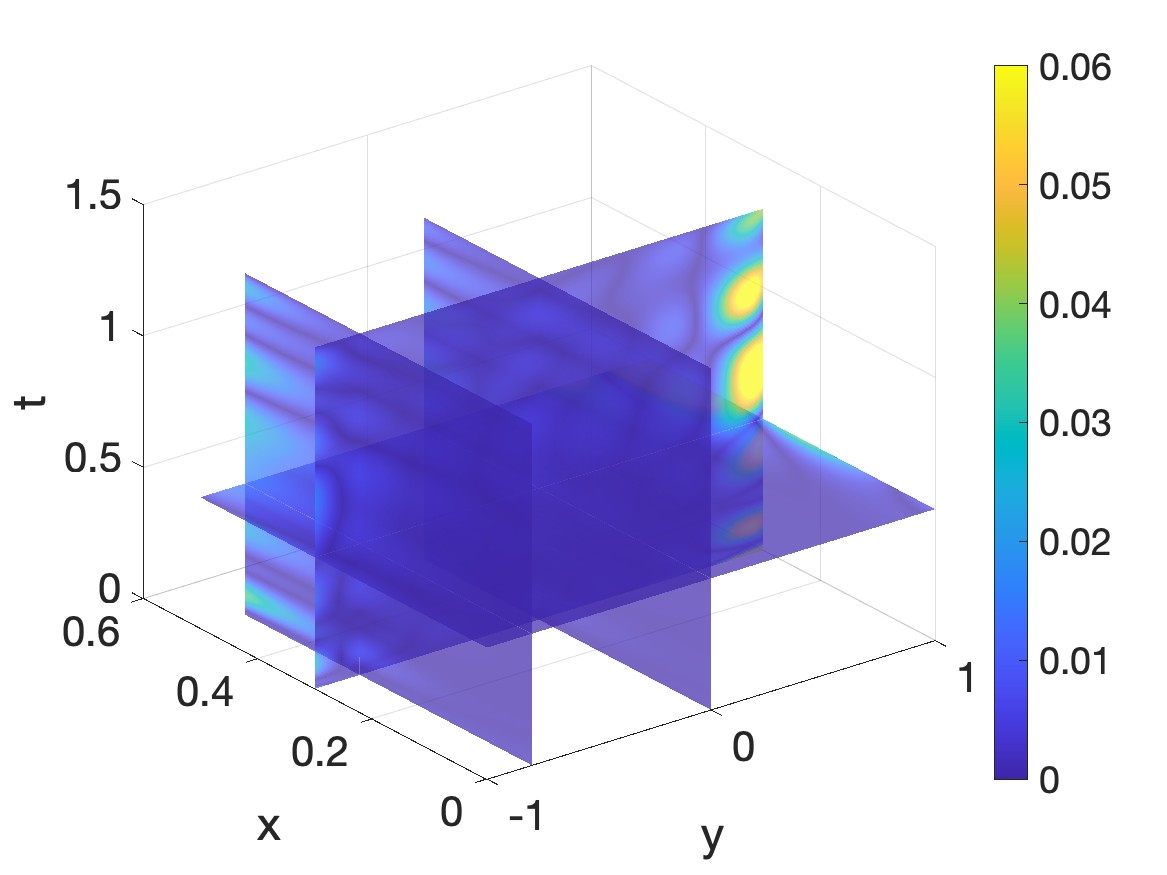}}
	\caption{\label{figtest3}Numerical results of Test 3 when the noise levels are $0,$ $5\%$, and $10\%$.
	The computation error primarily arises at the limit $x=0.5$.}
\end{figure}

Figure \ref{figtest3} provides visual evidence that Algorithm \ref{alg} generates acceptable solutions to \eqref{5.7}-\eqref{5.8}.
This unexpected performance is achieved although 
the nonlinearity in this test grows at the quadratic rate. 
Additionally, it is noteworthy that the computational error remains minimal, with the exception of areas near the boundary of the computational domain.
The relative errors $\frac{\|u_{\rm true} - u_{\rm comp}\|_{L^2}}{\|u_{\rm comp}\|_{L^2}}$ on $(0, 0.5) \times (-1, 1) \times (0, 1.3)$ are $0.016\%$, $2.11\%$, and $4.47\%$ when the noise levels are $\delta = 0\%,$ $\delta = 5\%,$ and $\delta = 10\%$, respectively.

\subsection{Test 4}

We solve a more complicated nonlinear parabolic equation
\begin{equation}
	u_t(\x) = \Delta u(\x, t) + |\nabla u(\x)|^2 + 5u(\x, t) 
	- \cos^2(2x + y + t)
	+ \cos(2x + y + t)
	\label{5.10}
\end{equation}
for all $\x = (x, y) \in \Omega$ and $t \in (0, T).$
The data on $\Gamma_T$ are given by
\begin{equation}
	g(y, t) = \sin(y + t) 
	\quad\mbox{and}
	\quad
	q(y, t) = 2 \cos(y + t)
	\label{5.11}
\end{equation}
for $y \in (-1, 1)$ and $t \in (0, 1.5).$

\begin{figure}[h!]
	\subfloat[The function $u_{\rm true}$]{\includegraphics[width = .3\textwidth]{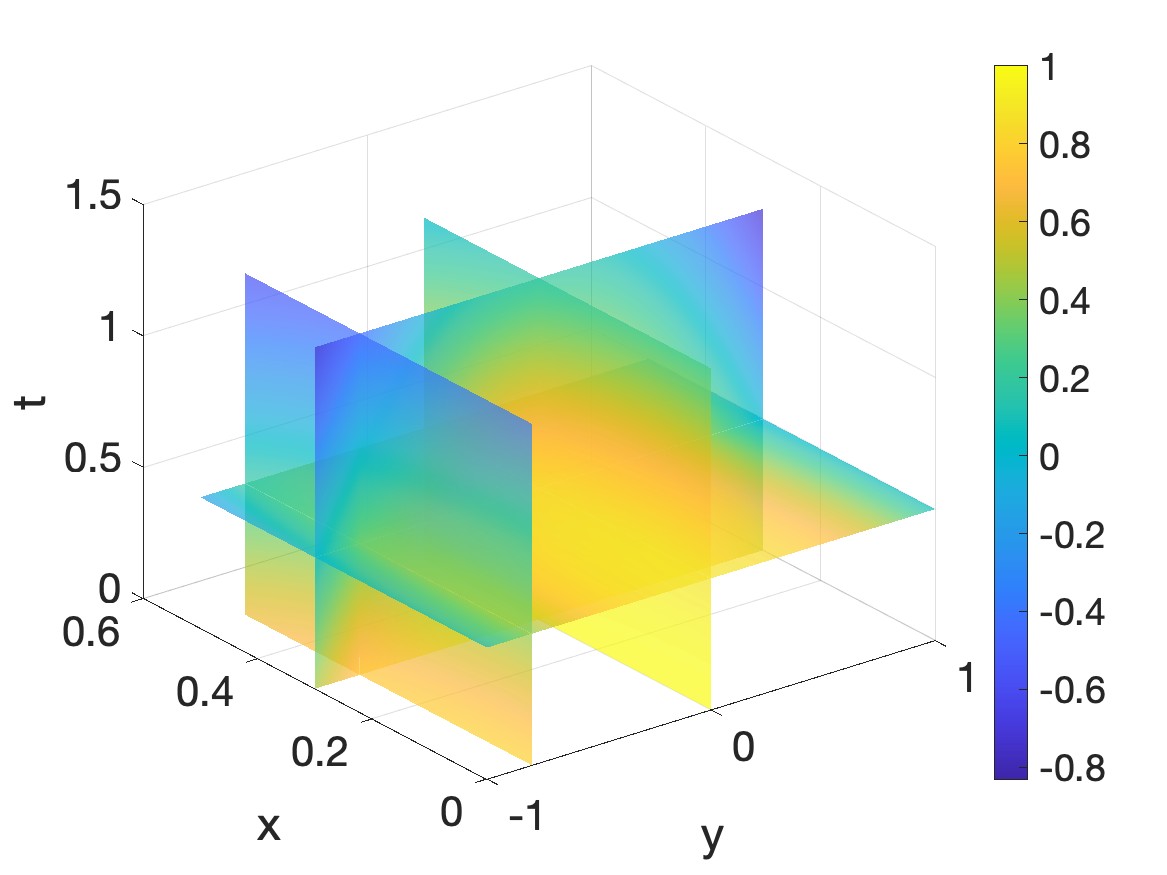}}
	
	\subfloat[The function $u_{\rm comp}$, $\delta = 0\%$]{\includegraphics[width = .3\textwidth]{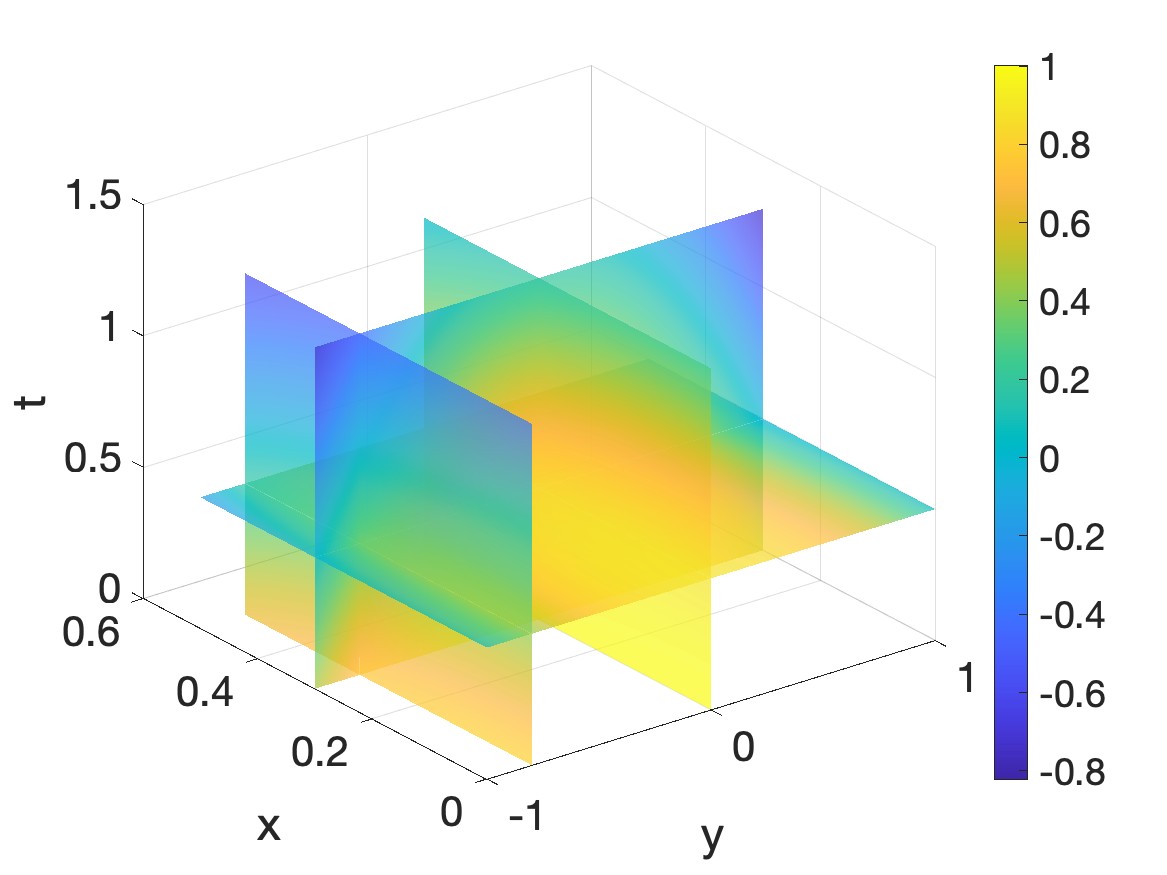}}
	\quad
	\subfloat[The function $u_{\rm comp}$, $\delta = 5\%$]{\includegraphics[width = .3\textwidth]{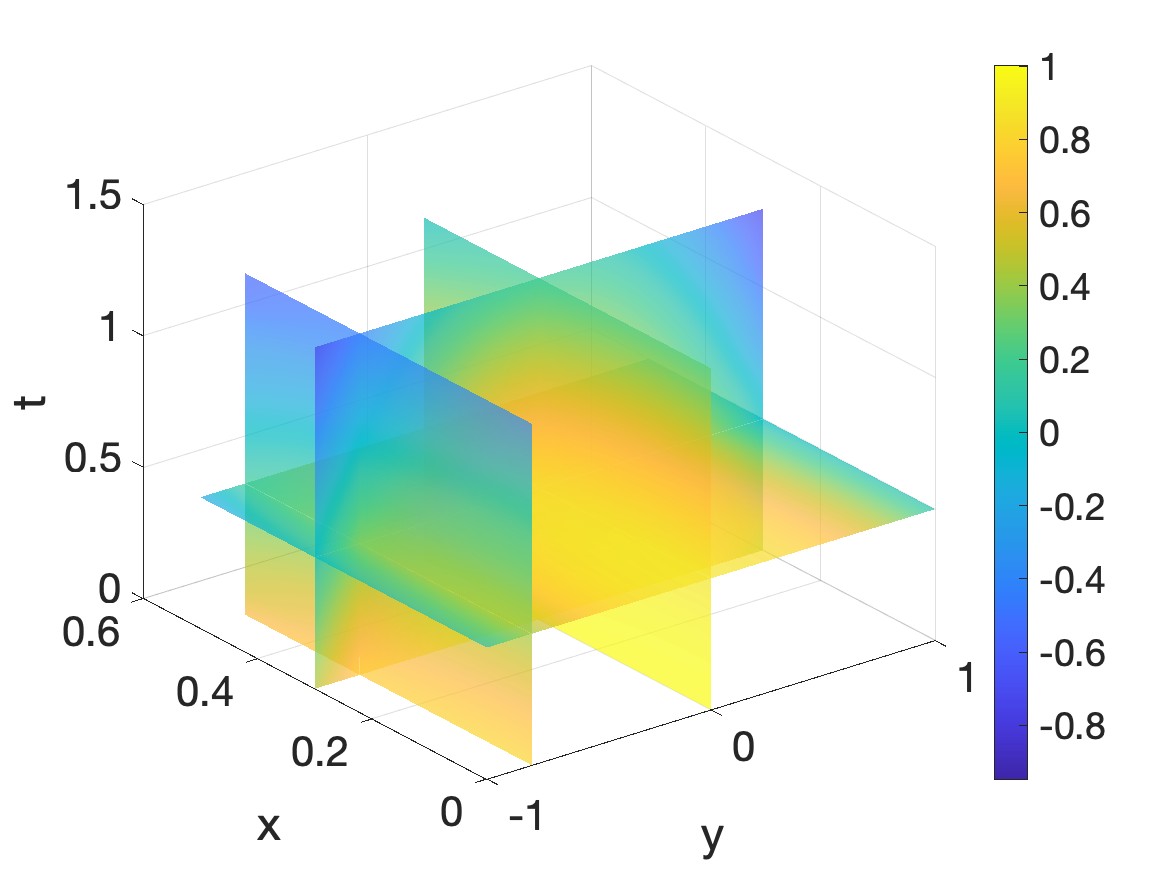}}
	\quad
	\subfloat[The function $u_{\rm comp}$, $\delta = 10\%$]{\includegraphics[width = .3\textwidth]{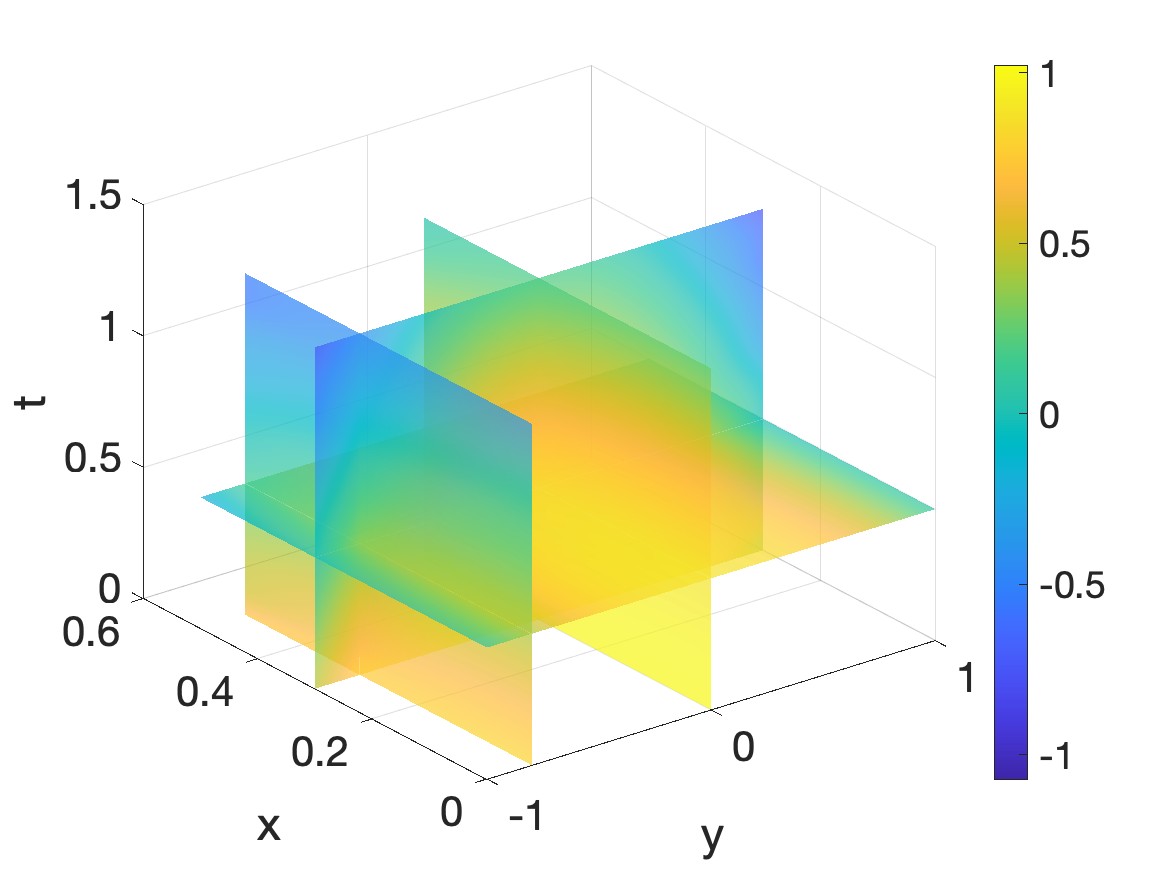}}
	
	\subfloat[The function $\frac{|u_{\rm true} - u_{\rm comp}|}{\|u_{\rm true}\|_{L^\infty}}$, $\delta = 0\%$]{\includegraphics[width = .3\textwidth]{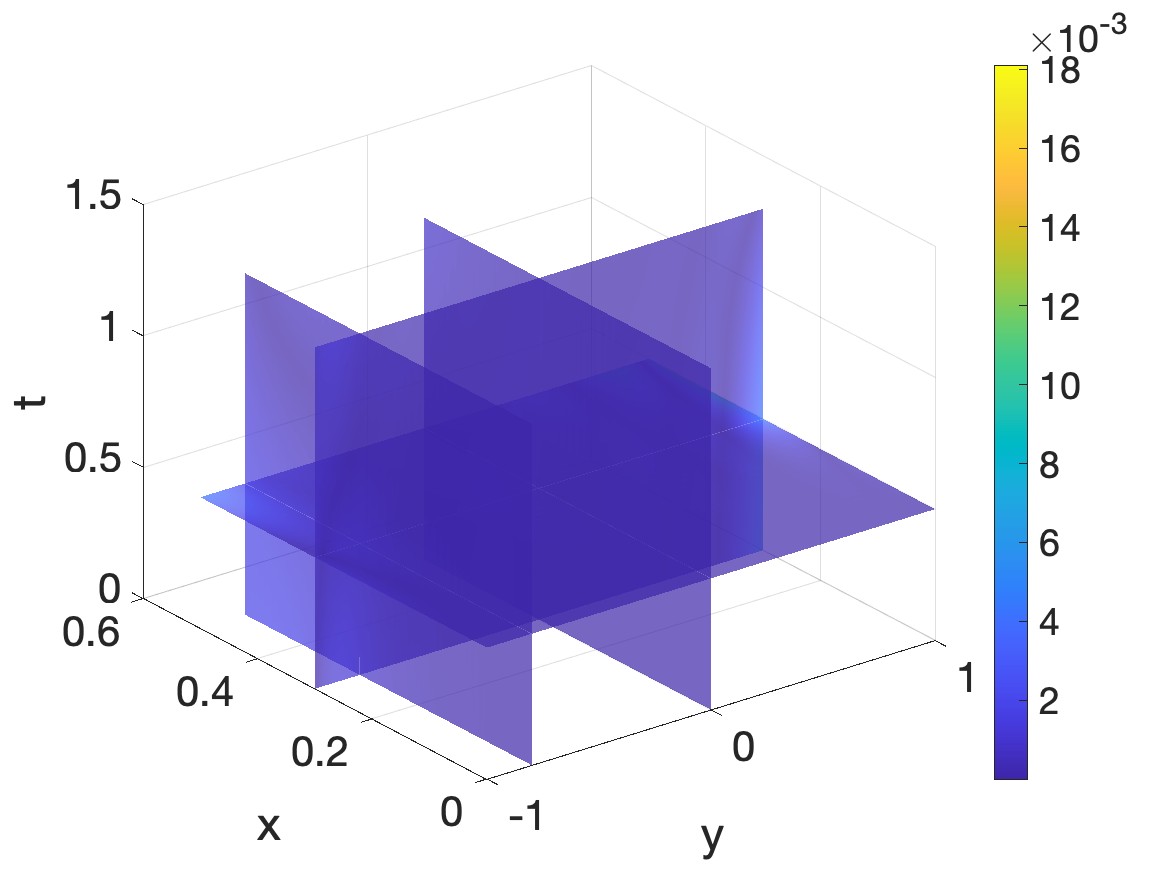}}
	\quad
	\subfloat[The function $\frac{|u_{\rm true} - u_{\rm comp}|}{\|u_{\rm true}\|_{L^\infty}}$, $\delta = 5\%$]{\includegraphics[width = .3\textwidth]{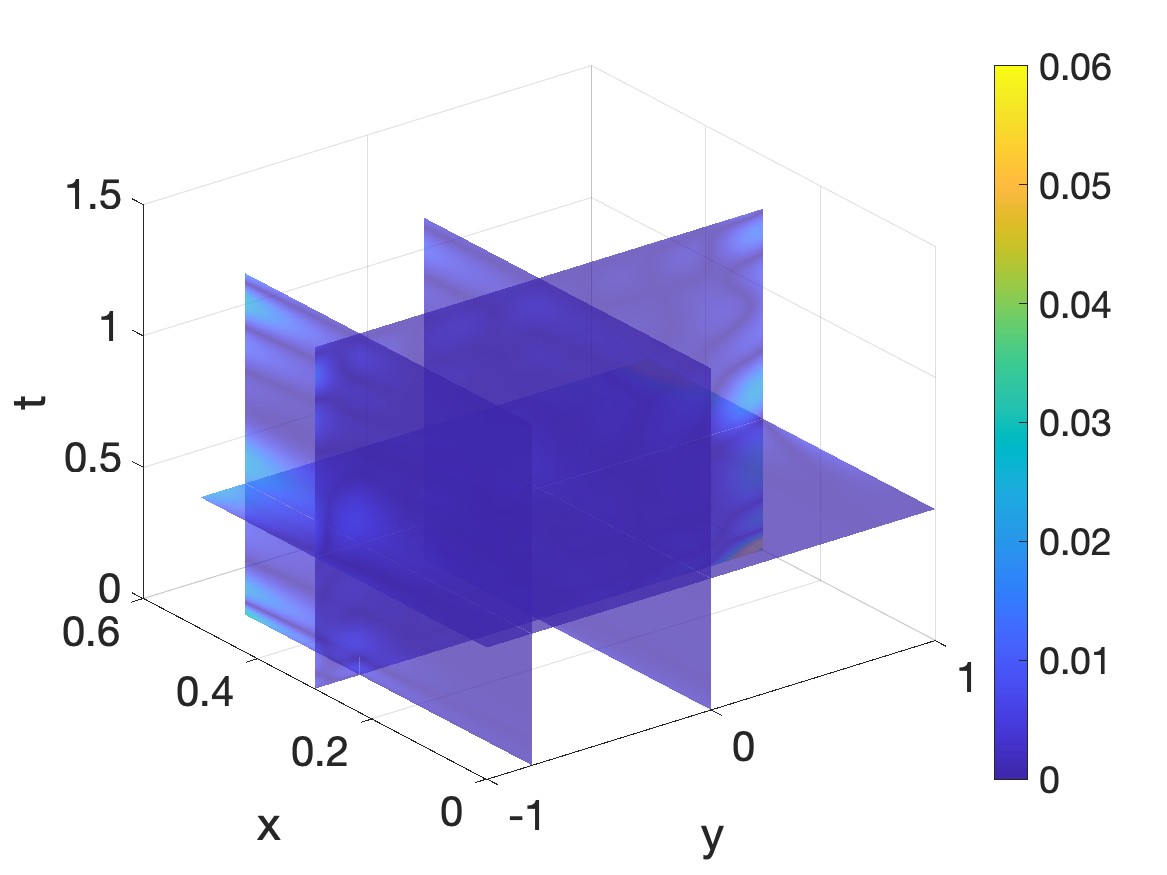}}
	\quad
	\subfloat[The function $\frac{|u_{\rm true} - u_{\rm comp}|}{\|u_{\rm true}\|_{L^\infty}}$, $\delta = 10\%$]{\includegraphics[width = .3\textwidth]{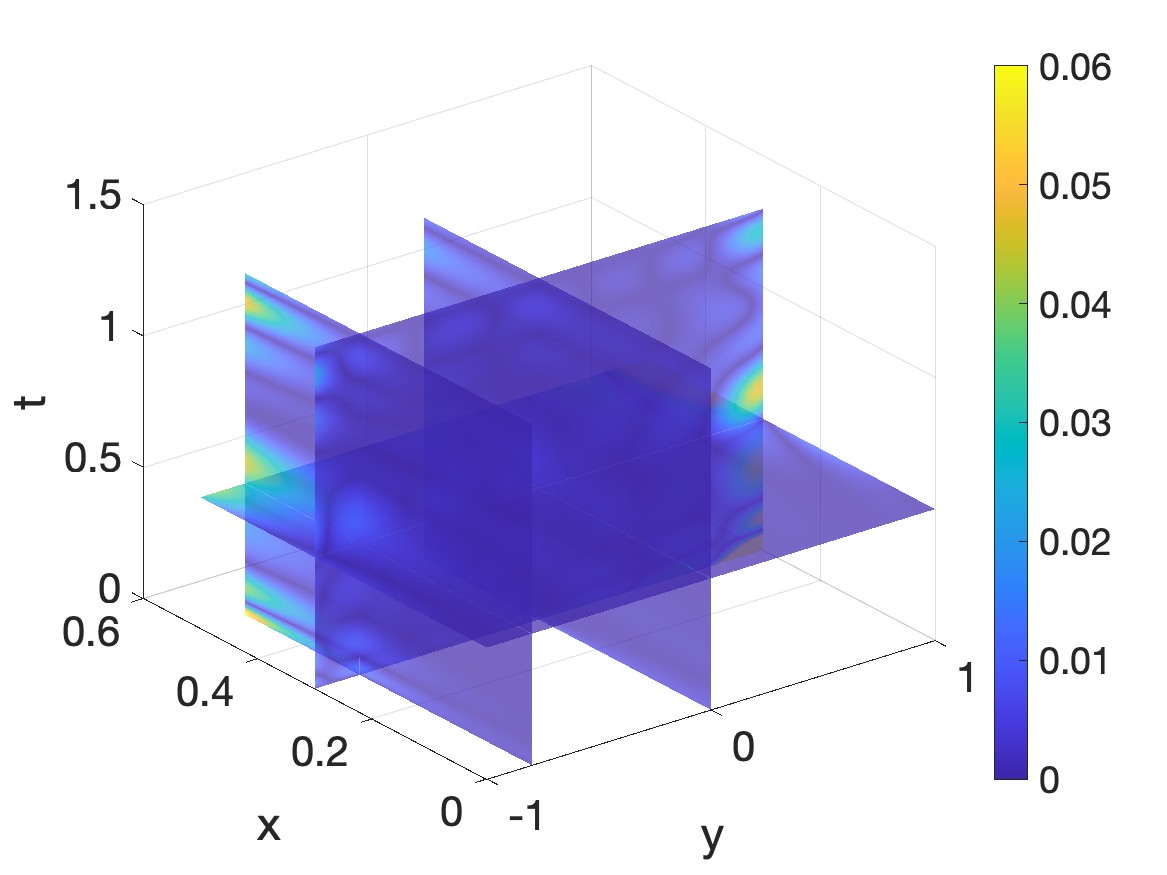}}
	\caption{\label{figtest4}Numerical results of Test 4 when the noise levels are $0,$ $5\%$, and $10\%$.
	The computation error primarily arises at the boundary of the computational domain $x \approx 0.5$, $y \approx \pm 1,$ $t \approx 1.5$.}
\end{figure}

Like in test 3, Figure \ref{figtest4} demonstrates that Algorithm \ref{alg} successfully produces satisfactory solutions for \eqref{5.10}-\eqref{5.11}, despite the nonlinearity growing at a quadratic rate. It is also worth mentioning that the computational error remains low, except in regions close to the computational domain's boundary. The relative errors $\frac{|u_{\rm true} - u_{\rm comp}|{L^2}}{|u{\rm comp}|_{L^2}}$ on the domain $(0, 0.5) \times (-1, 1) \times (0, 1.3)$ amount to $0.18\%$, $1.33\%$, and $2.62\%$ for noise levels of $\delta = 0\%,$ $\delta = 5\%,$ and $\delta = \%$, respectively.

\begin{Remark}
In addition to its computational accuracy, the dimensional reduction method possesses a significant advantage in terms of processing speed. We used a MacBook Pro laptop with a 1.4 GHz Quad-Core Intel Core i5 processor to compute the numerical results above. Algorithm \ref{alg} was able to deliver real-time solutions for Test 1 and Test 2. Moreover, it took only 12.5 seconds to compute solutions for Test 3 and Test 4. The remarkable speed is noteworthy as we perform numerical calculations for $(d + 1)$-dimensional solutions to nonlinear parabolic problems with minuscule increments in both spatial and temporal domains.
\end{Remark} 

\section{Concluding remarks} 
\label{sec remarks}

We have introduced a novel technique, known as the dimensional reduction method, for computing numerical solutions to nonlinear parabolic equations when lateral data is provided on one side, $\Gamma$, of the boundary of the computational domain $\Omega$. To tackle the nonlinearity and the lack of the data on $\partial \Omega \setminus \Gamma$, we employ a polynomial-exponential basis for the Fourier expansion of the solution.
By truncating the Fourier series, we can generate an approximate model in the form of a system of 1D ordinary differential equations. This system can be solved using the well-established Runge-Kutta method.
Impressively, this approach rapidly yields accurate numerical solutions.

 \section*{Acknowledgement}
The work by DNH was partially supported by VAST project NVCC01.06/22-23. The works of TTL and LHN were partially supported  by National Science Foundation grant DMS-2208159, and  by funds provided by the Faculty Research Grant program at UNC Charlotte Fund No. 111272,  and by
the CLAS small grant provided by the College of Liberal Arts \& Sciences, UNC Charlotte. 
The authors extend their gratitude to VIASM and the organizers of the PDE and Related Topics workshop held at VIASM in Hanoi, Vietnam, in July 2022. The initiation of this project took place at this meeting.


\begin{thebibliography}{10}

\bibitem{AbhishekLeNguyenKhan}
A.~Abhishek, T.~T. Le, L.~H. Nguyen, and T.~Khan.
\newblock The {C}arleman-{N}ewton method to globally reconstruct a source term
  for nonlinear parabolic equation.
\newblock {\em preprint, arXiv:2209.08011}, 2022.

\bibitem{Nehad1998}
N.~Al-Khalidy.
\newblock A general space marching algorithm for the solution of
  two-dimensional boundary inverse heat conduction problems.
\newblock {\em Numerical Heat Transfer, Part B: Fundamentals}, 34(3):339--360,
  1998.

\bibitem{alifanov-book1}
O.M. Alifanov.
\newblock {\em Inverse Heat Transfer Problems}.
\newblock Springer Verlag, Heidelberg, 1995.

\bibitem{alifanov-book2}
O.M. Alifanov, E.~A. Artiukhin, and S.~V. Rumiantsev.
\newblock {\em Extreme Methods for Solving Ill-Posed Problems With Applications
  to Inverse Heat Transfer Problems}.
\newblock Begell House, New York - Wallingford, 1995.

\bibitem{Arora:2019}
S.~Arora and J.~Dabas.
\newblock Inverse heat conduction problem in two-dimensional anisotropic
  medium.
\newblock {\em Int. J. Appl. Comput. Math}, 5(6):161, 2019.

\bibitem{beck-book}
J.V. Beck, B.~Blackwell, and C.R.St. Clair~Jr.
\newblock {\em Inverse Heat Conduction: Ill-Posed Problems}.
\newblock Wiley-Interscience publication. Wiley, 1985.

\bibitem{Burggraf1964}
O.~R. Burggraf.
\newblock An exact solution of the inverse problem in heat conduction theory
  and applications.
\newblock {\em Journal of Heat Transfer}, 86(3):373--380, 08 1964.

\bibitem{BadiaDuong:jiip2002}
A.~El~Badia and T.~Ha-Duong.
\newblock On an inverse source problem for the heat equation. application to a
  pollution detection problem.
\newblock {\em Journal of Inverse and Ill-posed Problems}, 10:585--599, 2002.

\bibitem{Felde2014}
I.~Felde and W.~Shi.
\newblock Hybrid approach for solution of inverse heat conduction problems.
\newblock In {\em 2014 IEEE International Conference on Systems, Man, and
  Cybernetics (SMC)}, pages 3896--3899, 2014.

\bibitem{Fisher:ae1937}
R.~A. Fisher.
\newblock The wave of advance of advantageous genes.
\newblock {\em Annals of Eugenics}, 7(4):355--369, 1937.

\bibitem{Ginsberg1963}
F.~Ginsberg.
\newblock On the {C}auchy problem for the one-dimensional heat equation.
\newblock {\em Mathematics of Computation}, 17:257--269, 1963.

\bibitem{Hao1992}
D.~N. H\`ao.
\newblock A noncharacteristic {C}auchy problem for linear parabolic equations
  {II}: a variational method.
\newblock {\em Numer. Funct. Anal. and Optimiz}, 13:541--564, 1992.

\bibitem{Hao1994}
D.~N. H\`ao.
\newblock A non-characteristic {C}auchy problem for linear parabolic equations
  and related inverse problems: {I}. {S}olvability.
\newblock {\em Inverse Problems}, 10:295--315, 1994.

\bibitem{haob}
D.~N. H\`ao.
\newblock {\em Methods for Inverse Heat Conduction Problems Problems}.
\newblock Peter Lang Verlag, Frankfurt/Main, Bern, New York, Paris, 1998.

\bibitem{isakov}
V.~Isakov.
\newblock {\em Inverse Problems in Partial Differential Equations}.
\newblock Springer, New York, Berlin, 1998.

\bibitem{JARNY19912911}
Y.~Jarny, M.N. Ozisik, and J.P. Bardon.
\newblock A general optimization method using adjoint equation for solving
  multidimensional inverse heat conduction.
\newblock {\em International Journal of Heat and Mass Transfer},
  34(11):2911--2919, 1991.

\bibitem{Klibanov:ip2006}
M.~V. Klibanov.
\newblock Estimates of initial conditions of parabolic equations and
  inequalities via lateral {C}auchy data.
\newblock {\em Inverse Problems}, 22:495--514, 2006.

\bibitem{Klibanov:jiip2017}
M.~V. Klibanov.
\newblock Convexification of restricted {D}irichlet to {N}eumann map.
\newblock {\em J. Inverse and Ill-Posed Problems}, 25(5):669--685, 2017.

\bibitem{KlibanovLiBook}
M.~V. Klibanov and J.~Li.
\newblock {\em Inverse Problems and Carleman Estimates: Global Uniqueness,
  Global Convergence and Experimental Data}.
\newblock De Gruyter, 2021.

\bibitem{KlibanovNguyenTran:JCP2022}
M.~V. Klibanov, L.~H. Nguyen, and H.~V. Tran.
\newblock Numerical viscosity solutions to {H}amilton-{J}acobi equations via a
  {C}arleman estimate and the convexification method.
\newblock {\em Journal of Computational Physics}, 451:110828, 2022.

\bibitem{LadyZhenskaya:ams1968}
O.~A. Ladyzhenskaya, V.A. Solonnikov, and N.~N. Ural'tseva.
\newblock {\em Linear and quasilinear equations of Parabolic Type}, volume~23.
\newblock American Mathematical Society, Providence, RI, 1968.

\bibitem{LeCON2023}
T.~T. Le.
\newblock Global reconstruction of initial conditions of nonlinear parabolic
  equations via the {C}arleman-contraction method.
\newblock In D-L. Nguyen, L.~H. Nguyen, and T-P. Nguyen, editors, {\em Advances
  in Inverse problems for Partial Differential Equations}, volume 784 of {\em
  Contemporary Mathematics}, pages 23--42. American Mathematical Society, 2023.

\bibitem{LeNguyen:jiip2022}
T.~T. Le and L.~H. Nguyen.
\newblock A convergent numerical method to recover the initial condition of
  nonlinear parabolic equations from lateral {C}auchy data.
\newblock {\em Journal of Inverse and Ill-posed Problems,}, 30(2):265--286,
  2022.

\bibitem{LeNguyen:JSC2022}
T.~T. Le and L.~H. Nguyen.
\newblock The gradient descent method for the convexification to solve boundary
  value problems of quasi-linear {PDEs} and a coefficient inverse problem.
\newblock {\em Journal of Scientific Computing}, 91(3):74, 2022.

\bibitem{LeNguyenNguyenPowell:JOSC2021}
T.~T. Le, L.~H. Nguyen, T-P. Nguyen, and W.~Powell.
\newblock The quasi-reversibility method to numerically solve an inverse source
  problem for hyperbolic equations.
\newblock {\em Journal of Scientific Computing}, 87:90, 2021.

\bibitem{LeNguyenTran:CAMWA2022}
T.~T. Le, L.~H. Nguyen, and H.~V. Tran.
\newblock A {C}arleman-based numerical method for quasilinear elliptic
  equations with over-determined boundary data and applications.
\newblock {\em Computers and Mathematics with Applications}, 125:13--24, 2022.

\bibitem{LiYamamotoZou:cpaa2009}
J.~Li, M.~Yamamoto, and J.~Zou.
\newblock Conditional stability and numerical reconstruction of initial
  temperature.
\newblock {\em Communications on Pure and Applied Analysis}, 8:361--382, 2009.

\bibitem{Loulou:2006}
T.~Loulou and E.~Scott.
\newblock An inverse heat conduction problem with heat flux measurements.
\newblock {\em International Journal for Numerical Methods in Engineering},
  67:1587--1616, 2006.

\bibitem{Mohebbi2021}
F.~Mohebbi, B.~Evans, and T.~Rabczuk.
\newblock Solving direct and inverse heat conduction problems in functionally
  graded materials using an accurate and robust numerical method.
\newblock {\em International Journal of Thermal Sciences}, 159:106629, 2021.

\bibitem{NguyenKlibanov:ip2022}
L.~H. Nguyen and M.~V. Klibanov.
\newblock Carleman estimates and the contraction principle for an inverse
  source problem for nonlinear hyperbolic equations.
\newblock {\em Inverse Problems}, 38:035009, 2022.

\bibitem{NguyenLeNguyenKlibanov:arxiv2023}
P.~M. Nguyen, T.~T. Le, L.~H. Nguyen, and M.~V. Klibanov.
\newblock Numerical differentiation by the polynomial-exponential basis.
\newblock {\em preprint arXiv:2304.05909}, 2023.

\bibitem{Ozisik2017}
M.~N. \"Ozi\c{s}ik, H.~R.~B. Orlande, M.~J. Cola\c{c}o, and R.~M. Cotta.
\newblock {\em Finite Difference Methods in Heat Transfer}.
\newblock CRC Press, 2nd edition, 2017.

\bibitem{Roy2023}
A.~D. Roy and S.~K. Dhiman.
\newblock Solutions of one-dimensional inverse heat conduction problems: a
  review.
\newblock {\em Transactions of the Canadian Society for Mechanical Engineering,
  https://doi.org/10.1139/tcsme-2022-0157}, 2023.

\end{thebibliography}
\end{document}